\documentclass[preprint,12pt]{elsarticle}

\usepackage{graphicx}

\usepackage{amssymb}

\usepackage[figuresright]{rotating}
\usepackage{subfigure}
\usepackage{url}

\newtheorem{theorem}{Theorem}

\newtheorem{definition}{Definition}
\newtheorem{example}{Example}

\usepackage{amssymb}
\usepackage{amsmath}
\usepackage{amsfonts}
\usepackage{epstopdf}
\usepackage{listings}
\usepackage{mathptmx}
\usepackage[T1]{fontenc}
\usepackage[scaled]{helvet}
\usepackage{setspace}

\usepackage{graphicx}
\usepackage{xcolor}
\usepackage[nottoc,notlof,notlot]{tocbibind}
\usepackage{booktabs}
\usepackage{float}
\restylefloat{table}

\usepackage[pagebackref,linktocpage,breaklinks,colorlinks,%
                linkcolor=black,anchorcolor=black,citecolor=black,%
                filecolor=black,menucolor=black,runcolor=black,%
                urlcolor=black]{hyperref}

\journal{to a journal}

\begin{document}

\begin{frontmatter}



\title{Locally Linearized Runge Kutta method of Dormand and Prince}


\author[JCJ,JCJ1]{J.C. Jimenez}
\author[AS]{A. Sotolongo}
\author[JMSB]{J.M. Sanchez-Bornot}

\address[JCJ]{Instituto de Cibern\'{e}tica, Matem\'{a}tica y
F\'{\i}sica, La Habana, Cuba. e-mail: jcarlos@icimaf.cu}
\address[JCJ1]{The Institute of Statistical Mathematics, Tokyo, Japan}
\address[AS]{Universidad de la Habana, La
Habana, Cuba. e-mail: alina.sotolongo@yahoo.com}
\address[JMSB]{Centro de
Neurociencias de Cuba, La Habana, Cuba. e-mail: bornot@gmail.com}

\begin{abstract}
In this paper, the effect that produces the local linearization of
the embedded Runge-Kutta formulas of Dormand and Prince for initial
value problems is studied. For this, embedded Locally Linearized
Runge-Kutta formulas are defined and their performance is analyzed
by means of exhaustive numerical simulations. For a variety of
well-known physical equations with different dynamics, the
simulation results show that the locally linearized formulas exhibit
significant higher accuracy than the original ones, which implies a
substantial reduction of the number of time steps and, consequently,
a sensitive reduction of the overall computation cost of their
adaptive implementation.
\end{abstract}


\end{frontmatter}


\section{Introduction}
\label{Introduction} It is well known (see, e.g., \cite{Cartwright
1992, Stewart 1992, Skufca 2004}) that conventional numerical
schemes such as Runge-Kutta, Adams-Bashforth, predictor-corrector
and others produce misleading dynamics when integrating Ordinary
Differential Equations (ODEs). Typical problems are, for instance,
the convergence to spurious steady states, changes in the basis of
attraction, appearance of spurious bifurcations, etc. This might
yield serious mistakes in the interpretation and analysis of the
processes under consideration in practical control engineering or in
applied sciences. The essence of such difficulties is that the
dynamic of the numerical schemes (considered as discrete dynamical
systems) is far richer than that of its continuous counterparts.
Contrary to the popular belief, drawbacks of this type may no be
solved by reducing the stepsize of the numerical method. Therefore,
it is highly desirable the development of numerical integrators that
preserve, as much as possible, the dynamical properties of the
underlaying dynamical system for all step sizes or relative big
ones. In this direction, some modest advances has been achieved by a
number of relative recent integrators of the class of Exponential
Methods, which are characterized by the explicit use of exponentials
to obtain an approximate solution. An example of such integrators
are the High Order Local Linearization (HOLL) methods based on
Runge-Kutta schemes \cite{de la Cruz 06, de la Cruz 11, Jimenez 09}.

HOLL integrators are obtained by splitting, at each time step, the
solution of the original ODE in two parts: the solution of a linear
ODE plus the solution of an auxiliary ODE. The linear equation is
solved by a Local Linearization (LL) scheme \cite{Jimenez02 AMC,
Jimenez05 AMC} in such a way that A-stability is ensured, whereas
the solution of the auxiliary one can be approximated by any
conventional numerical integrator, preferably a high order explicit
scheme. Originally, HOLL methods were introduced as a flexible
approach for increasing the order of convergence of the order-$2$ LL
method but, in addition, they can be thought as a strategy for
constructing high order A-stable explicit schemes based on
conventional explicit integrators. For this reason, if we focus on
the conventional integrator involved in a particular HOLL scheme,
then it is natural to say that the first one has been locally
linearized. In this way, if a Runge Kutta scheme is used to
approximate the above mentioned auxiliary ODE, the resulting HOLL
scheme are indistinctly called Local Linearization - Runge Kutta
(LLRK) scheme or Locally Linearized Runge Kutta (LLRK) scheme.

In \cite{de la Cruz 11}, general results on the convergence,
stability and dynamical properties of the Locally Linearized Runge
Kutta method were studied. Specifically, it was demonstrated that:
1) the LLRK approach defines a general class of high order A-stable
explicit integrators that preserve the convergence rate of the
involved (not A-stable) explicit RK schemes; 2) in contrast with
others A-stable explicit methods (such as Rosenbrock or the
Exponential integrators), the RK coefficients involved in the LLRK
integrators are not constrained by any stability condition and they
just need to satisfy the usual, well-known order conditions of RK
schemes, which makes the LLRK approach more flexible and simple; 3)
LLRK integrators have a number of convenient dynamical properties
such as the linearization preserving and the conservation of the
exact solution dynamics around hyperbolic equilibrium points and
periodic orbits; and 4) because of the flexibility in the numerical
implementation of the LLRK discretizations, specific-purpose LLRK
schemes can be designed for certain classes of ODEs, e.g., for
moderate or large systems of equations. On the other hand,
simulation studies carried out in \cite{de la Cruz 06, de la Cruz
11, Sotolongo 11} have shown that, for a variety of test equations,
LLRK schemes of order $3$ and $4$ preserve much better the stability
and dynamical properties of the actual solutions than their
corresponding conventional RK schemes.

However, the accuracy and computational efficiency of the Local
Linearization methods have been much less considered up to now,
being the dynamical properties of such schemes the focus of previous
studies and the main reason for the development of these methods.
The few available results are the following. On an identical time
partition \cite{de la Cruz 11}, the LLRK scheme based on the
classical order-$4$ RK scheme displays better accuracy than the
order-$5$ RK formula of Dormand \& Prince \cite{Dormand80} in the
integration of a variety of ODEs. On different time partitions \cite
{Sotolongo 11}, similar results are obtained by an adaptive
implementation of the mentioned LLRK scheme in comparison with the
Matlab code ode$45$, which provides an adaptive implementation of
the embedded RK formulas of Dormand \& Prince. However, this is
achieved at expense of additional evaluations of the vector field,
and with larger overall computational time. With this respect, the
main drawback of that adaptive LLRK scheme is the absence of a
computationally efficient strategy based on embedded formulas.

The main purpose of this work is introducing an adaptive LLRK scheme
based on the embedded RK formulas of Dormand \& Prince and
evaluating, with simulations, its accuracy and computational
efficiency in order to study the effect that the local linearization
produces on these known RK formulas. The Matlab code developed with
this goal is, same as the Matlab code ode$45$, addressed to low
dimensional non stiff initial value problems for medium to low
accuracies.

The paper is organized as follows. In the Section 2, a basic
introduction on the Local Linearization - Runge Kutta (LLRK) schemes
is presented. In the Section 3, the embedded Locally Linearized
Runge-Kutta formulas are defined, and an adaptive implementation of
them is described. In the last two sections, the results of a
variety of exhaustive numerical simulations with well-known test
equations are presented and discussed respectively.

\section{Notations and preliminaries}

Let $\mathcal{D}\subset \mathbb{R}^{d}$ be an open set. Consider the $d$%
-dimensional differential equation%
\begin{eqnarray}
\frac{d\mathbf{x}\left( t\right) }{dt} &=&\mathbf{f}\left( t,\mathbf{x}%
\left( t\right) \right) \text{, \ \ }t\in \left[ t_{0},T\right]
\label{ODE-LLA-1} \\
\mathbf{x}(t_{0}) &=&\mathbf{x}_{0},  \label{ODE-LLA-2}
\end{eqnarray}%
where $\mathbf{x}_{0}\in \mathcal{D}$ is a given initial point, and $\mathbf{%
f}:\left[ t_{0},T\right] \times \mathcal{D}\longrightarrow
\mathbb{R}
^{d}$ is a differentiable function. Lipschitz and smoothness
conditions on the function $\mathbf{f}$ are assumed in order to
ensure a unique solution of this equation in $\mathcal{D}$.

Let $\left( t\right) _{h}=\left\{ t_{n}:n=0,1,\ldots ,N\right\} $ be
a time discretization with maximum stepsize $h$ defined as a
sequence of times that
satisfy the conditions $t_{0}<t_{1}<\cdots <t_{N}=T$ and $%
\sup\limits_{n}(h_{n})\leq h<1$, where $h_{n}=t_{n+1}-t_{n}$ for
$n=0,\ldots ,N-1$.

For a given $(t_{n},\mathbf{y}_{n})$, let $\mathbf{v}_{n+1}=\mathbf{y}_{n}+%
\Lambda_{1}\left( t_{n},\mathbf{y}_{n};h_{n}\right) $ be an order-$%
\gamma _{1}$ approximation to solution of the linear ODE
\begin{eqnarray}
\frac{d\mathbf{z}_{1}\left( t\right) }{dt} &=&\mathbf{B}_{n}\mathbf{z}%
_{1}(t)+\mathbf{b}_{n}\left( t\right) \mathbf{,}\text{ \ \ }t\in
[t_{n},t_{n+1}],  \label{ODE-LLS-13} \\
\mathbf{z}_{1}\left( t_{n}\right) &=&\mathbf{y}_{n}
\label{ODE-LLS-13b}
\end{eqnarray}%
at $t_{n+1}$, and let $\mathbf{w}_{n+1}=\Lambda%
_{2}^{\mathbf{y}_{n}}\left( t_{n},\mathbf{0};h_{n}\right) $ be an order-$%
\gamma _{2}$ Runge-Kutta scheme approximating the solution of the
nonlinear ODE
\begin{eqnarray}
\frac{d\mathbf{z}_{2}\left( t\right) }{dt}
&=&\mathbf{q}(t_{n},\mathbf{y}
_{n};t\mathbf{,\mathbf{z}}_{2}\left( t\right) \mathbf{),}%
\text{ \ \ }t\in [t_{n},t_{n+1}],\quad  \label{ODE-LLS-14} \\
\mathbf{z}_{2}\left( t_{n}\right) &=&\mathbf{0}  \label{ODE-LLS-14b}
\end{eqnarray}%
at $t_{n+1}$, where $\mathbf{B}_{n}=\mathbf{f}_{\mathbf{x}}\left( t_{n},%
\mathbf{y}_{n}\right) $ is a $d\times d$ constant matrix, and
\begin{equation*}
\mathbf{b}_{n}(t)=\mathbf{f}_{t}\left( t_{n},\mathbf{y}_{n}\right) (t-t_{n})+%
\mathbf{f}\left( t_{n},\mathbf{y}_{n}\right)
-\mathbf{B}_{n}\mathbf{y}_{n}
\end{equation*}%
and
\begin{eqnarray*}
\mathbf{q}(t_{n},\mathbf{y}_{n};s\mathbf{,\xi} )& = & \mathbf{f}(s,\mathbf{y}_{n}+\Lambda_{1}\left( t_{n},\mathbf{y}_{n};s-t_{n}\right) +\mathbf{\xi })-%
\mathbf{f}_{\mathbf{x}}(t_{n},\mathbf{y}_{n})\Lambda_{1}\left(
t_{n},\mathbf{y}_{n};s-t_{n}\right)  \\
& & -\mathbf{f}_{t}\left( t_{n},\mathbf{y}%
_{n}\right) (s-t_{n})-\mathbf{f}\left( t_{n},\mathbf{y}_{n}\right)
\end{eqnarray*}
are $d$-dimensional vectors. Here, $\mathbf{f}_{x}$ and
$\mathbf{f}_{t}$ denote the partial derivatives respect to
$\mathbf{x}$ and $t$, respectively. Note that the vector field of
the equation (\ref{ODE-LLS-14}) not only depends on the point
$(t_{n},\mathbf{y}_{n})$ but also of the numerical flow used to
approximate $\mathbf{z}_{1}(t)$.

\begin{definition}
\label{Definition LLRK schemes}(\cite{de la Cruz 06, Jimenez 09, de
la Cruz
11}) A Local Linearization - Runge Kutta (LLRK) scheme for the ODE (\ref%
{ODE-LLA-1})-(\ref{ODE-LLA-2}) is defined by the recursive
expression
\begin{equation}
\mathbf{y}_{n+1}=\mathbf{y}_{n}+\Lambda_{1}\left( t_{n},\mathbf{y}%
_{n};h_{n}\right) +\Lambda_{2}^{\mathbf{y}_{n}}\left( t_{n},%
\mathbf{0};h_{n}\right)  \label{LLRK scheme}
\end{equation}%
for all $t_{n}\in $ $\left( t\right) _{h}$, starting with $\mathbf{y}_{0}=%
\mathbf{x}_{0}$.
\end{definition}

Local truncation error, rate of convergence and various dynamical
properties of the general class of Local Linearization - Runge Kutta
schemes (\ref{LLRK scheme}) has been studied in \cite{de la Cruz
11}.

According to the Definition \ref{Definition LLRK schemes}, a variety
of LLRK schemes can be derived. In previous works \cite{de la Cruz
06, de la Cruz 11} , the Local Linearization scheme based on
Pad\'{e} approximations \cite {Jimenez02 AMC, Jimenez05 AMC} has
been used to integrate the linear ODE
(\ref{ODE-LLS-13})-(\ref{ODE-LLS-13b}), whereas the so called
\textit{four order classical }Runge-Kutta scheme \cite{Butcher 2008}
has been applied to integrate the nonlinear ODE
(\ref{ODE-LLS-14})-(\ref{ODE-LLS-14b}). This yields the order-$4$
LLRK scheme
\begin{equation}
\mathbf{y}_{n+1}=\mathbf{y}_{n}+\mathbf{u}_{4}+\frac{h_{n}}{6}(2\mathbf{k}%
_{2}+2\mathbf{k}_{3}+\mathbf{k}_{4}),  \label{LLRK4 scheme}
\end{equation}%
where%
\begin{equation*}
\mathbf{u}_{j}=\mathbf{L}(\mathbf{P}_{6,6}(2^{-\kappa _{j}}\mathbf{D}%
_{n}c_{j}h_{n}))^{2^{\kappa _{j}}}\mathbf{r}
\end{equation*}%
and\ \ \
\begin{equation*}
\mathbf{k}_{j}=\mathbf{f}\left( t_{n}+c_{j}h_{n},\mathbf{y}_{n}+\mathbf{u}%
_{j}+c_{j}h_{n}\mathbf{k}_{j-1}\right) -\mathbf{f}\left( t_{n},\mathbf{y}%
_{n}\right) -\mathbf{f}_{\mathbf{x}}\left(
t_{n},\mathbf{y}_{n}\right) \mathbf{u}_{j}\ -\mathbf{f}_{t}\left(
t_{n},\mathbf{y}_{n}\right) c_{j}h_{n},
\end{equation*}%
with $\mathbf{k}_{1}\equiv \mathbf{0}$ and $c=\left[
\begin{array}{cccc}
0 & \frac{1}{2} & \frac{1}{2} & 1%
\end{array}%
\right] $. Here, $\mathbf{P}_{p,q}(\cdot )$ denotes the
$(p,q)$-Pad\'{e} approximation for exponential matrices \cite{Moler
1978}, and $\kappa _{j}$
the smallest integer number such that $\left\Vert 2^{-\kappa _{j}}\mathbf{D}%
_{n}c_{j}h_{n}\right\Vert \leq \frac{1}{2}$. The matrices $\mathbf{D}_{n}$, $%
\mathbf{L}$ and $\mathbf{r}$ are defined as
\begin{equation*}
\mathbf{D}_{n}=\left[
\begin{array}{ccc}
\mathbf{f}_{\mathbf{x}}(t_{n},\mathbf{y}_{n}) & \mathbf{f}_{t}(t_{n},\mathbf{%
y}_{n}) & \mathbf{f}(t_{n},\mathbf{y}_{n}) \\
0 & 0 & 1 \\
0 & 0 & 0%
\end{array}%
\right] \in \mathbb{R}^{(d+2)\times (d+2)},
\end{equation*}%
$\mathbf{L}=\left[
\begin{array}{ll}
\mathbf{I}_{d} & \mathbf{0}_{d\times 2}%
\end{array}%
\right] $ and $\mathbf{r}^{\intercal }=\left[
\begin{array}{ll}
\mathbf{0}_{1\times (d+1)} & 1%
\end{array}%
\right] $ for non-autonomous ODEs; and as
\begin{equation*}
\mathbf{D}_{n}=\left[
\begin{array}{cc}
\mathbf{f}_{\mathbf{x}}(\mathbf{y}_{n}) & \mathbf{f}(\mathbf{y}_{n}) \\
0 & 0%
\end{array}%
\right] \in \mathbb{R}^{(d+1)\times (d+1)},
\end{equation*}%
$\mathbf{L}=\left[
\begin{array}{ll}
\mathbf{I}_{d} & \mathbf{0}_{d\times 1}%
\end{array}%
\right] $ and $\mathbf{r}^{\intercal }=\left[
\begin{array}{ll}
\mathbf{0}_{1\times d} & 1%
\end{array}%
\right] $ for autonomous equations.

On an identical time partition \cite{de la Cruz 11}, LLRK formula
(\ref {LLRK4 scheme}) displays better accuracy than the order-$5$ RK
formula of Dormand \& Prince in the integration of a variety of
ODEs. On different time partitions \cite{Sotolongo 11}, similar
results are obtained by an adaptive implementation of LLRK formula
(\ref{LLRK4 scheme}) in comparison with the Matlab ode$45$ code,
which provides an adaptive implementation of the embedded RK
formulas of Dormand \& Prince. However, this is achieved at expense
of additional evaluations of the vector field $\mathbf{f}$, and with
larger overall computational time. This cost can be sensitively
reduced by using the\ $(2,2)$-Pad\'{e} approximations instead of the
order $(6,6)$ one used in formula (\ref{LLRK4 scheme}), preserving
the order of convergence and without significant lost of accuracy
\cite{Sotolongo 11}.

Local truncation error, rate of convergence, A-stability and various
dynamical properties of the LLRK schemes based on Pad\'{e}
approximations has also been studied in \cite{de la Cruz 11}.

For a precise comparison of the Local Linearization - Runge Kutta
method with well-known integration methods such as Rosenbrock,
Exponential Integrators, Splitting Methods and others, the
interested reader might read \cite{de la Cruz 11} or \cite{Jimenez
09}.

\section{Numerical scheme}

\subsection{Embedded Locally Linearized Runge-Kutta formulas}

In view of the Definition \ref{Definition LLRK schemes}, new
integration formulas can be obtained as follows. Similarly to the
LLRK scheme (\ref {LLRK4 scheme}), the Local Linearization scheme
based on Pad\'{e} approximations \cite{Jimenez02 AMC, Jimenez05 AMC}
is used for integrating the linear ODE
(\ref{ODE-LLS-13})-(\ref{ODE-LLS-13b}) but, instead of the classical
order-$4$ RK scheme, the embedded Runge-Kutta formulas of Dormand \&
Prince \cite{Dormand80} is now applied to integrate the nonlinear
ODE (\ref{ODE-LLS-14})-(\ref{ODE-LLS-14b}). This yields the embedded
Locally Linearized Runge-Kutta formulas

\begin{equation}
\mathbf{y}_{n+1}=\mathbf{y}_{n}+\mathbf{u}_{s}+h_{n}\sum_{j=1}^{s}b_{j}%
\mathbf{k}_{j}\ \ \ \text{and}\ \ \ \ \ \widehat{\mathbf{y}}_{n+1}=\mathbf{y}%
_{n}+\mathbf{u}_{s}+h_{n}\sum_{j=1}^{s}\widehat{b}_{j}\mathbf{k}_{j},
\label{embedded LLRK45}
\end{equation}%
where $s=7$ is the number of the stages,
\begin{equation}
\mathbf{u}_{j}=\mathbf{L}(\mathbf{P}_{p,q}(2^{-\kappa _{j}}\mathbf{D}%
_{n}c_{j}h_{n}))^{2^{\kappa _{j}}}\mathbf{r} \label{embedded
LLRK45b}
\end{equation}%
and\ \ \
\begin{equation*}
\mathbf{k}_{j}=\mathbf{f(}t_{n}+c_{j}h_{n},\mathbf{y}_{n}+\mathbf{u}%
_{j}+h_{n}\sum_{i=1}^{s-1}a_{j,i}\mathbf{k}_{i})-\mathbf{f}\left( t_{n},%
\mathbf{y}_{n}\right) -\mathbf{f}_{\mathbf{x}}\left( t_{n},\mathbf{y}%
_{n}\right) \mathbf{u}_{j}\ -\mathbf{f}_{t}\left( t_{n},\mathbf{y}%
_{n}\right) c_{j}h_{n},
\end{equation*}%
with $\mathbf{k}_{1}\equiv \mathbf{0}$ and Runge-Kutta coefficients $a_{j,i}$%
, $b_{j}$, $\widehat{b}_{j}$ and $c_{j}$ defined in the Table
\ref{Table EmbeddedRK45}. Here, $\mathbf{P}_{p,q}(\cdot )$ denotes
the $(p,q)$-Pad\'{e} approximation for exponential matrices with
$p+q>4$. The number $\kappa _{j}$ and the matrices $\mathbf{D}_{n}$,
$\mathbf{L}$ and $\mathbf{r}$ are defined as in the previous
section.

\begin{table}[h]
  \centering
\begin{tabular}{p{1.5cm}@{\vrule height 15pt depth 8pt width 0pt}|p{1.5cm}@{\vrule height 15pt depth 8pt width 0pt}p{1.5cm}@{\vrule height 15pt depth 8pt width 0pt}p{1.5cm}@{\vrule height 15pt depth 8pt width 0pt}p{1.5cm}@{\vrule height 15pt depth 8pt width 0pt}p{1.5cm}@{\vrule height 15pt depth 8pt width 0pt}p{1.5cm}@{\vrule height 15pt depth 8pt width 0pt}p{1.5cm}@{\vrule height 15pt depth 8pt width 0pt}p{1.5cm}@{\vrule height 15pt depth 8pt width 0pt}p{1.5cm}@{\vrule height 15pt depth 8pt width 0pt}}
  0 &  &  &  &  &  &  &   \\
  $\frac{1}{5}$ & $\frac{1}{5}$ &  &  &  &  &  &   \\
  $\frac{3}{10}$ & $\frac{3}{40}$ & $\frac{9}{40}$ &  &  &  &  &   \\
  $\frac{4}{5}$ & $\frac{44}{45}$ & $-\frac{56}{15}$ & $\frac{32}{9}$ &  &  &  &  \\
  $\frac{8}{9}$ & $\frac{19372}{6561}$ & $-\frac{25360}{2187}$ & $\frac{64448}{6561}$ & $-\frac{212}{729}$ &  &  &  \\
  $1$ & $\frac{9017}{3168}$ & $-\frac{355}{33}$ & $\frac{46732}{5247}$ & $\frac{49}{176}$ & $-\frac{5103}{18656}$ &  &  \\
  $1$ & $\frac{35}{384}$ & $0$ & $\frac{500}{1113}$ & $\frac{125}{192}$ & $-\frac{2187}{6784}$ & $\frac{11}{84}$ &  \\ \hline
  $y$ & $\frac{35}{384}$ & $0$ & $\frac{500}{1113}$ & $\frac{125}{192}$ & $-\frac{2187}{6784}$ & $\frac{11}{84}$ & $0$ \\ \hline
  $\hat{y}$ & $\frac{5179}{57600}$ & $0$ & $\frac{7571}{16695}$ & $\frac{393}{640}$ & $-\frac{92097}{339200}$ & $\frac{187}{2100}$ & $\frac{1}{40}$  \\
  \hline
\end{tabular}
  \caption{Coefficients tableau for the embedded formulas.}\label{Table
EmbeddedRK45}
\end{table}

The local truncation error, the rate of convergence and the
A-stability of the LLRK formulas (\ref{embedded LLRK45}) will be
consider in what follows. With this purpose, these formulas are
rewritten as
\[
\mathbf{y}_{n+1}=\mathbf{y}_{n}+h_{n}\varphi (t_{n,}\mathbf{y}_{n};h_{n})%
\text{\ \ \ \ \ \ \ \ \ and \ \ \ \ \ \ \ \ \ \ \ }\widehat{\mathbf{y}}%
_{n+1}=\mathbf{y}_{n}+h_{n}\widehat{\varphi }(t_{n,}\mathbf{y}_{n};h_{n})%
\text{,}
\]%
and the following additional notations are introduced. Let
$\mathcal{D}$
be an open subset of $\mathbb{R}^{d}$, $\mathcal{M}$ an upper bound for $%
\left\Vert \mathbf{f}_{\mathbf{x}}\right\Vert $ on $[t_{0},T]\times \mathcal{%
D}$, and $\mathcal{L}$ the Lipschitz constant of the function $%
q(t,x(t);\cdot )$ (which exists for all $t\in \lbrack t_{0},T]$
because Lemma 6 in \cite{de la Cruz 11} under regular conditions for
$\mathbf{f}$). Denote by $L_{n+1}$ the local truncation
error of the Local Linearization scheme $\mathbf{y}_{n+1}=\mathbf{y}%
_{n}+u_{4}$ when it is applied to the linear equation (\ref{ODE-LLS-13})-(%
\ref{ODE-LLS-13b}), for which the inequality
\[
L_{n+1}\leq Ch^{p+q+1}
\]%
holds with positive constant $C$ \cite{Jimenez05 AMC}. Further, denote by $%
L_{n+1}^{1}$ and $L_{n+1}^{2}$ the local truncation errors of the
classical embedded Runge-Kutta formulas of Dormand and Prince when
they are applied to the nonlinear equation
(\ref{ODE-LLS-14})-(\ref{ODE-LLS-14b}), for which the inequalities
\[
L_{n+1}^{1}\leq C_{1}h^{6}\text{\ \ \ \ \ \ \ \ \ and \ \ \ \ \ \ \ \ \ \ \ }%
L_{n+1}^{2}\leq C_{2}h^{5}
\]%
hold with positive constants $C_{1}$ and $C_{2}$
\cite{Dormand80,Hairer-Wanner93}.

\begin{theorem}\label{Main Theorem}
Let $\mathbf{x}$ be the solution of the ODE (\ref{ODE-LLA-1})-(\ref%
{ODE-LLA-2}) with vector field $\mathbf{f}$ six times continuously
differentiable on $[t_{0},T]\times \mathcal{D}$. Then, the embedded
Locally Linearized Runge-Kutta formulas (\ref{embedded LLRK45}) have
local truncation errors
\[
\left\Vert \mathbf{x}(t_{n+1})-\mathbf{x}(t_{n})-h_{n}\varphi (t_{n,}\mathbf{%
x}(t_{n});h_{n})\right\Vert \leq Kh_{n}^{p+q+1}+C_{1}h_{n}^{6}
\]%
and
\[
\left\Vert \mathbf{x}(t_{n+1})-\mathbf{x}(t_{n})-h_{n}\widehat{\varphi }%
(t_{n,}\mathbf{x}(t_{n});h_{n})\right\Vert \leq
Kh_{n}^{p+q+1}+C_{2}h_{n}^{5};
\]%
and global errors
\[
\left\Vert \mathbf{x}(t_{n+1})-\mathbf{y}_{n+1}\right\Vert \leq
M_{1}h^{\min \{p+q,5\}}
\]%
and
\[
\left\Vert \mathbf{x}(t_{n+1})-\widehat{\mathbf{y}}_{n+1}\right\Vert
\leq M_{2}h^{\min \{p+q,4\}}
\]%
for all $t_{n+1}\in \left( t\right) _{h}$ and $h$ small enough, where $K=C(1+%
\frac{\mathcal{M}}{\mathcal{L}}(e^{\mathcal{L}}-1))$ is a positive
constant, and $M_{1}$ and $M_{2}$ as well. In addition, the embedded
Locally Linearized Runge-Kutta formulas (\ref{embedded LLRK45}) are
A-stable if in the involved $(p,q)$-Pad\'{e} approximation the
inequality $p\leq q\leq p+2$ holds$.$
\end{theorem}

\textbf{Proof.} The local truncation errors and the global errors
are a straightforward consequence of Theorem 15 in \cite{de la Cruz
11}, whereas the A-stability is a direct result of Theorem 17 in
\cite{de la Cruz 11}. $\square $

Clearly, according to this result, the Locally Linearized
Runge-Kutta formulas (\ref{embedded LLRK45}) preserve the
convergence rate of the classical embedded Runge-Kutta formulas of
Dormand and Prince if the inequality $p+q>4$ holds. Further, note
that these Locally Linearized formulas not only preserve the
stability of the linear ODEs when $p\leq q\leq p+2$, but they are
also able to "exactly" (up to the precision of the floating-point
arithmetic) integrate this class of equations when $p+q=12$ (for the
numerical precision of the current personal computers \cite{Moler
1978}).

In addition, and trivially, the embedded Locally Linearized
Runge-Kutta formulas (\ref{embedded LLRK45}) inherit the dynamical
properties derived in \cite{de la Cruz 11} for the general class of
Local Linearization - Runge Kutta methods.

\subsection{Adaptive strategy\label{Adpative scheme}}

In order to write a code that automatically adjust the stepsizes for
achieving a prescribed tolerance of the local error at each step, an
adequate adaptive strategy is necessary. At glance, the automatic
stepsize control for the embedded RK formulas of Dormand \& Prince
seems to fit well for the embedded LLRK formulas (\ref{embedded
LLRK45}). In what follows, the adaptive strategy of the Matlab code
ode$45$ for these formulas is described.

Once the values for the relative and absolute tolerances $RTol$ and
$ATol$, and for the maximum and minimum stepsizes $h_{\max }$ and
$h_{\min }$ are set, the basic steps of the algorithm are:

\begin{enumerate}
\item Estimation of the initial stepsize
\begin{equation*}
h_{0} = \min \{h_{\max },\max \{h_{\min },\Delta \}\}
\end{equation*}
where
\begin{equation*}
\Delta =\left\{
\begin{array}{cc}
\frac{1}{r_{h}} & \text{if }h_{max}\cdot r_{h}>1 \\
h_{\max } & \text{otherwise}%
\end{array}%
\right.
\end{equation*}%
with%
\begin{equation*}
r_{h}=\frac{1}{0.8\cdot RTol^{1/5}}\max_{i=1...d}\left\{ \frac{\mathbf{f}%
^{i}(\mathbf{y}_{0})}{\max \left\{ |\mathbf{y}_{0}^{i}|,tr\right\}
}\right\}
\end{equation*}%
and $tr=$ $\frac{ATol}{RTol}$. Initialize $fail=0$.

\item Evaluation of the embedded formula (\ref{embedded LLRK45})

\item Estimation of the error
\begin{equation*}
error={\LARGE ||}\frac{\mathbf{y}_{n+1}-\widehat{\mathbf{y}}_{n+1}}{\underset%
{i=1,\ldots ,d}{\max }\left\{ |\mathbf{y}_{n}^{i}|,|\ \mathbf{y}%
_{n+1}^{i}|,tr\right\} }{\LARGE ||}_{\infty }
\end{equation*}

\item Estimation of a new stepsize%
\begin{equation*}
h_{new}=\min \{h_{\max },\max \{h_{\min },\Delta \}\}
\end{equation*}%
where
\begin{equation*}
\Delta =\left\{
\begin{array}{cc}
0.8\cdot {\large (}\frac{RTol}{error}{\large )}^{1/5}\cdot h & \text{if }%
error\leq RTol \\
\max \{0.1,0.8\cdot {\large (\frac{RTol}{error})}^{1/5}\}\cdot h & \text{if }%
error>RTol\text{ and }fail=0 \\
0.5\cdot h & \text{if }error>RTol\text{ and }fail=1%
\end{array}%
\right.
\end{equation*}

\item Validation of $y_{n+1}$: if $error\leq RTol$, then accept $y_{n+1}$ as
an approximation to $\mathbf{x}$ at $t_{n+1}=t_{n}+h$. Otherwise,
return to 2 with $h_{n}=h_{new}$ and $fail=1.$

\item Control of the final step: if $t_{n}+h=T$, stop. If $t_{n}+h+h_{new}>T$%
, then redefine $h_{new}=T-(t_{n}+h)$.

\item Return to 2 with $n=n+1$, $h_{n}=h_{new}$, and $fail=0$.
\end{enumerate}

\subsection{Continuous formula}

Continuous formulas of RK methods are usually defined for computing
the solutions on a dense set of time instants with minimum
computational cost. Typically \cite{Hairer-Wanner93}, they are
constructing by means of a
polynomial interpolation of the RK formulas between two consecutive times $%
t_{n},t_{n+1}\in $ $\left( t\right) _{h}$.

By a simple combination of the LLRK formulas (\ref{embedded LLRK45})
with
the continuous formulas of the Dormand \& Prince RK method \cite%
{Hairer-Wanner93} for (\ref{ODE-LLS-14})-(\ref{ODE-LLS-14b}), a
continuous $7$-stage LLRK formula can be defined as

\begin{equation}
\mathbf{y}(t_{n}+\theta h_{n})=\mathbf{y}_{n}+\mathbf{u}(\theta
h_{n})+h_{n}\sum_{j=1}^{7}b_{j}(t_{n}+\theta h_{n})\mathbf{k}_{j}\
,\ \ 0<\theta <1,  \label{continuousLLRK$45$}
\end{equation}%
for all $t_{n},t_{n+1}\in $ $\left( t\right) _{h}$, where
\begin{equation}
\mathbf{u}(\theta h_{n})=\mathbf{L}(\mathbf{P}_{p,q}(2^{-\kappa
_{j}}\mathbf{D}_{n}\theta h_{n}))^{2^{\kappa _{j}}}\mathbf{r}
\label{continuousLLRK$45$b}
\end{equation}%
is a $d$-dimensional vector, and
\begin{equation*}
b_{j}(\delta )=\sum\limits_{i=1}^{4}\alpha _{i,j}\delta ^{i}
\end{equation*}%
is a polynomial with coefficients $\alpha _{i,j}$. Here, the function $%
\mathbf{k}_{j}$, the matrices $\mathbf{D}_{n}$, $\mathbf{L}$ and $\mathbf{r}$%
, and the number $\kappa _{j}$ are defined as in (\ref{embedded
LLRK45}), as well as the $(p,q)$-Pad\'{e} approximation
$\mathbf{P}_{p,q}$. The coefficients $\alpha _{i,j}$, defined in
Table \ref{Table continuous RK}, coincide with those of the
continuous RK formula implemented in the Matlab code ode$45$.

\begin{table}[tbp] \centering
$%
\begin{tabular}{|c|c|c|c|c|}
\hline $j/i$ & $1$ & $2$ & $3$ & $4$ \\ \hline $1$ & $1$ & $-183/64$
& $37/12$ & $-145/128$ \\ \hline $2$ & $0$ & $0$ & $0$ & $0$ \\
\hline $3$ & $0$ & $1500/371$ & $-1000/159$ & $1000/371$ \\ \hline
$4$ & $0$ & $-125/32$ & $125/12$ & $-375/64$ \\ \hline $5$ & $0$ &
$9477/3392$ & $-729/106$ & $25515/6784$ \\ \hline $6$ & $0$ &
$-11/7$ & $11/3$ & $-55/28$ \\ \hline $7$ & $0$ & $3/2$ & $-4$ &
$5/2$ \\ \hline
\end{tabular}%
$%
\caption{Values of the coefficient $\alpha _{i,j}$ involved in the
continuous LLRK formula (\ref{continuousLLRK$45$}) \label{Table
continuous RK}.}
\end{table}%

\subsection{LLDP45 code\label{LLRK45 codes}}

This subsection describes a Matlab2007b(32bits) implementation of
the adaptive scheme described above, which will be denoted as
LLDP$45$ code.

In order to make a fair comparison between the linearized and the
nonlinearized RK formulas, the LLDP$45$ code is an exact copy of the
ode$45$ one with the exception of the program lines corresponding to
the embedded and continuous formulas of Dormand and
Prince, which are replaced by the formulas (\ref{embedded LLRK45}) and (\ref%
{continuousLLRK$45$}) respectively. We recall that the code ode$45$
implements the adaptive strategy of the subsection \ref{Adpative
scheme}\ for the embedded RK formulas of Dormand \& Prince, which is
considered by many authors the most recommendable code to be applied
as first try for most problems \cite{Shampine97}.

Note that, the embedded LLRK formulas (\ref{embedded LLRK45})
require the computation of six Pad\'{e} approximations
$\mathbf{P}_{p,q}$ at each integration step, which increases the
computational cost of the original embedded RK formulas.
Nevertheless, this number of Pad\'{e} approximations
can be reduced by taking in to account that: a) $(\mathbf{P}%
_{p,q}(2^{-\kappa }\mathbf{D}_{n}c_{j}h_{n}))^{2^{\kappa }}$ gives
an approximation to exponential matrix
$\mathbf{e}^{\mathbf{D}_{n}c_{j}h_{n}}$; and b) the flow property of
the exponential operator. Indeed, this can be carried out in two
steps:

\begin{enumerate}
\item approximating $\mathbf{e}^{\mathbf{D}_{n}h_{n}/90}$ by the matrix $%
\mathbf{M}_{1/90}=$ $(\mathbf{P}_{p,q}(2^{-\kappa }\mathbf{D}%
_{n}h_{n}/90))^{2^{\kappa }},$ where $\kappa $ is the smallest
integer number such that $\left\Vert 2^{-\kappa
}\mathbf{D}_{n}h_{n}/90\right\Vert \leq \frac{1}{2}$; and

\item the successive computation of the matrices
\begin{eqnarray*}
\mathbf{M}_{2/90} &=&\mathbf{M}_{1/90}\mathbf{M}_{1/90}\ \ \ \ \ \ \
\ \
\mathbf{M}_{4/90}=\mathbf{M}_{2/90}\mathbf{M}_{2/90} \\
\mathbf{M}_{8/90} &=&\mathbf{M}_{4/90}\mathbf{M}_{4/90}\ \ \ \ \ \ \
\
\mathbf{M}_{16/90}=\mathbf{M}_{8/90}\mathbf{M}_{8/90} \\
\mathbf{M}_{32/90} &=&\mathbf{M}_{16/90}\mathbf{M}_{16/90}\ \ \ \ \
\
\mathbf{M}_{80/90}=\mathbf{M}_{32/90}\mathbf{M}_{16/90}\mathbf{M}_{32/90} \\
\mathbf{M}_{1/10} &=&\mathbf{M}_{8/90}\mathbf{M}_{1/90}\ \ \ \ \ \ \
\ \ \ \
\mathbf{M}_{1/5}=\mathbf{M}_{1/10}\mathbf{M}_{1/10} \\
\mathbf{M}_{2/5} &=&\mathbf{M}_{1/5}\mathbf{M}_{1/5}\ \ \ \ \ \ \ \
\ \ \ \
\ \ \mathbf{M}_{4/5}=\mathbf{M}_{2/5}\mathbf{M}_{2/5} \\
\mathbf{M}_{3/10} &=&\mathbf{M}_{1/10}\mathbf{M}_{1/5}\ \ \ \ \ \ \
\ \ \ \ \ \ \ \ \mathbf{M}_{1}=\mathbf{M}_{4/5}\mathbf{M}_{1/5}.
\end{eqnarray*}
\end{enumerate}

Consequently, the matrix $\mathbf{M}_{c_{j}}$ corresponding to each
RK
coefficient $c_{j}$\ provides an approximation to $\mathbf{e}^{\mathbf{D}%
_{n}c_{j}h_{n}}$, for all $j=1,..,6$. In this way, at each
integration step,
the code LLDP$45$ performs six evaluation of $\mathbf{f}$ (same than the ode$%
45$ code), one Jacobian matrix and one matrix exponential.

The matrix $\mathbf{M}_{1/90}$ is computed by means of the function "\textit{%
expmf}", which provides a C++ implementation of the classical $(p,q)$-Pad%
\'{e} approximations algorithm for exponential matrices with scaling
and squaring strategy \cite{Moler 1978}, and $p=q=3$.

\section{Numerical simulations}

In this section, the performance of the LLDP$45$ and ode$45$ codes
is compared by means of numerical simulations. To do so, a variety
of ODEs and
simulation types were selected. For all of them, the Relative Error%
\begin{equation}
RE=\underset{{\small i=1,\ldots ,d;}\text{ }{\small t}_{j}{\small \in (t)}%
_{h}}{\max }\left\vert \frac{\mathbf{x}^{i}(t_{j})-\mathbf{y}^{i}(t_{j})}{%
\mathbf{x}^{i}(t_{j})}\right\vert  \label{Relative Error}
\end{equation}%
between the "exact" solution $\mathbf{x}$ and its approximation
$\mathbf{y}$ is evaluated.

The simulations with the code ode$45$ were carried out with a wide
range of tolerances: crude with $RTol=10^{-3}$ and $ATol=10^{-6}$,
mild with $RTol=10^{-6}$ and $ATol=10^{-9}$, and refined with
$RTol=10^{-9}$ and $ATol=10^{-12}$. The Matlab code ode$15s$ with
refined tolerance $RTol=10^{-13}$ and $ATol=10^{-13}$ was used to
compute the "exact" solution $\mathbf{x}$ in all simulations.

\subsection{Test examples}

The first four examples have the semi-lineal form%
\begin{equation}
\frac{d\mathbf{x}}{dt}=\mathbf{Ax+g}(\mathbf{x),}
\label{Semi-Linear}
\end{equation}%
where $\mathbf{A}$ is a square matrix and $\mathbf{g}$ is a function of $%
\mathbf{x}$. The vector field of the first two examples have
Jacobian with eigenvalues on or near the imaginary axis, which made
these equations difficult to be integrated by conventional schemes
\cite{Shampine97}. The other two are also hard for conventional
explicit schemes since they are examples of stiff equations
\cite{Shampine97}. Example \ref{Example SNL} has an additional
complexity for a number of integrators that do not update the
Jacobian of the vector field at each integration step \cite%
{Shampine97,Hochbruck-etal09}: the Jacobian of the linear term has
positive eigenvalues, which results a problem for the integration in
a neighborhood of the stable equilibrium point $\mathbf{x}=1$.

\begin{example}
\label{Example PL} \ Periodic linear \cite{de la Cruz 11}%
\begin{equation*}
\frac{d\mathbf{x}}{dt}=\mathbf{A}(\mathbf{x}+2),
\end{equation*}%
with%
\begin{equation*}
\mathbf{A}=\left[
\begin{array}{cc}
i & 0 \\
0 & -i%
\end{array}%
\right] ,
\end{equation*}%
$\mathbf{x}(t_{0})=(-2.5,-1.5)$ and $[t_{0},T]=[0,4\pi ]$.
\end{example}

\begin{example}
\label{Example PNL}Periodic linear plus nonlinear term \cite{de la Cruz 11}%
\begin{equation*}
\frac{d\mathbf{x}}{dt}=\mathbf{A}(\mathbf{x}+2)+0.1\mathbf{x}^{2},
\end{equation*}%
where the matrix $\mathbf{A}$ is defined as in the previous example, $%
\mathbf{x}(t_{0})=(1,1)$, and $[t_{0},T]=[0,4\pi ]$.
\end{example}

\begin{example}
\label{Example SL}Stiff linear \cite{de la Cruz 11}$\ $%
\begin{equation*}
\frac{d\mathbf{x}}{dt}=-100\mathbf{H}(\mathbf{x+1}),
\end{equation*}%
where $\mathbf{H}$ is the 12-dimensional Hilbert matrix (with
conditioned number $1.69\times 10^{16}$), $\mathbf{x}^{i}(t_{0})=1$,
$i=1\ldots 12$, and $[t_{0},T]=[0,1].$
\end{example}

\begin{example}
\label{Example SNL}Stiff linear plus nonlinear term$\ $\cite{de la Cruz 11}%
\begin{equation*}
\frac{d\mathbf{x}}{dt}=100\mathbf{H}(\mathbf{x}-\mathbf{1})+100(\mathbf{x}-%
\mathbf{1})^{2}-60(\mathbf{x}^{3}-\mathbf{1}),
\end{equation*}%
where $\mathbf{H}$ is the 12-dimensional Hilbert matrix, $\mathbf{x}%
^{i}(t_{0})=-0.5$, $i=1\ldots 12$, and $[t_{0},T]=[0,1]$.
\end{example}

The following examples are well known nonlinear test equations. This
include highly oscillatory, non stiff and mild stiff equations.

\begin{example}
\label{FermiPastaUlam}Fermi--Pasta--Ulam equation defined by the
Hamiltonian
system \cite{Hairer-Lubich-Wanner06}%
\begin{equation*}
H(\mathbf{p,q})=\frac{1}{2}\sum\limits_{i=1}^{3}(\mathbf{p}_{2i-1}^{2}+%
\mathbf{p}_{2i}^{2})+\frac{w^{2}}{4}\sum\limits_{i=1}^{3}(\mathbf{q}_{2i}-%
\mathbf{q}_{2i-1})^{2}+\sum\limits_{i=0}^{3}(\mathbf{q}_{2i+1}-\mathbf{q}%
_{2i})^{2}
\end{equation*}%
with $w=50$, initial conditions $1,1,1/w,1$ for the four first
variables and zero for the remainder eight, and $[t_{0},T]=[0,15].$
\end{example}

\begin{example}
\label{Brusselator} Brusselator equation\ \cite{Hairer-Wanner93}:
\begin{eqnarray*}
\frac{dx_{1}}{dt} &=&1+{\normalsize x}_{1}^{2}x_{2}-4x_{1} \\
\frac{dx_{2}}{dt} &=&3x_{1}-x_{1}^{2}x_{2},
\end{eqnarray*}%
where $(x_{1}(t_{0}),x_{2}(t_{0}))=(1.5{\normalsize ,3)}$ and $%
[t_{0},T]=[0,20]$.
\end{example}

\begin{example}
\label{Rigid} Rigid body equation \cite{Hairer-Wanner93}:%
\begin{eqnarray*}
\frac{d{\normalsize x}_{1}}{dt} &=&{\normalsize x}_{2}x_{3} \\
\frac{d{\normalsize x}_{2}}{dt} &=&-{\normalsize x}_{1}x_{3} \\
\frac{d{\normalsize x}_{3}}{dt} &=&-0.51x_{1}{\normalsize x}_{2}
\end{eqnarray*}%
with $(x_{1}(t_{0}),x_{2}(t_{0}),x_{3}(t_{0}))=(0,1,1)$ over $%
[t_{0},T]=[0,12]$.\
\end{example}

\begin{example}
\label{Chemical reaction} Chemical reaction \cite{Shampine97}:%
\begin{eqnarray*}
\frac{d{\normalsize x}_{1}}{dt} &=&1.3(x_{3}-x_{1})+10400k({\normalsize x}%
_{1}){\normalsize x}_{2} \\
\frac{d{\normalsize x}_{2}}{dt} &=&1880(x_{4}-x_{2}(1+k({\normalsize x}_{1})%
\mathbf{))} \\
\frac{d{\normalsize x}_{3}}{dt} &=&1752-269x_{3}+267x_{1} \\
\frac{d{\normalsize x}_{4}}{dt} &=&0.1+320x_{2}-321x_{4}
\end{eqnarray*}%
where $k({\normalsize x_{1})=e^{(20.7-\frac{1500}{x_{1}})}}$. With
initial condition $(50,0,600,0.1)$ over $[t_{0},T]=[0,1]$, this is
mild stiff equation.
\end{example}

\begin{example}
\label{Van der Pol}Van der Pol equation \cite{Hairer-Wanner96}:
\begin{eqnarray*}
\frac{d{\normalsize x}_{1}}{dt} &=&x_{2} \\
\frac{d{\normalsize x}_{2}}{dt} &=&\varepsilon
(1-x_{2}^{2})x_{1}+x_{2}
\end{eqnarray*}%
with $(x_{1}(t_{0}),x_{2}(t_{0}))=(2,0)$. With $\varepsilon =1$ and $%
\varepsilon =10^{2}$, this is a non stiff and a mild stiff equation
on the intervals $[t_{0},T]=[0,20]$ and $[t_{0},T]=[0,300]$,
respectively.
\end{example}

\begin{figure}
\subfigure[PerLin]{\includegraphics[width=0.32\textwidth]{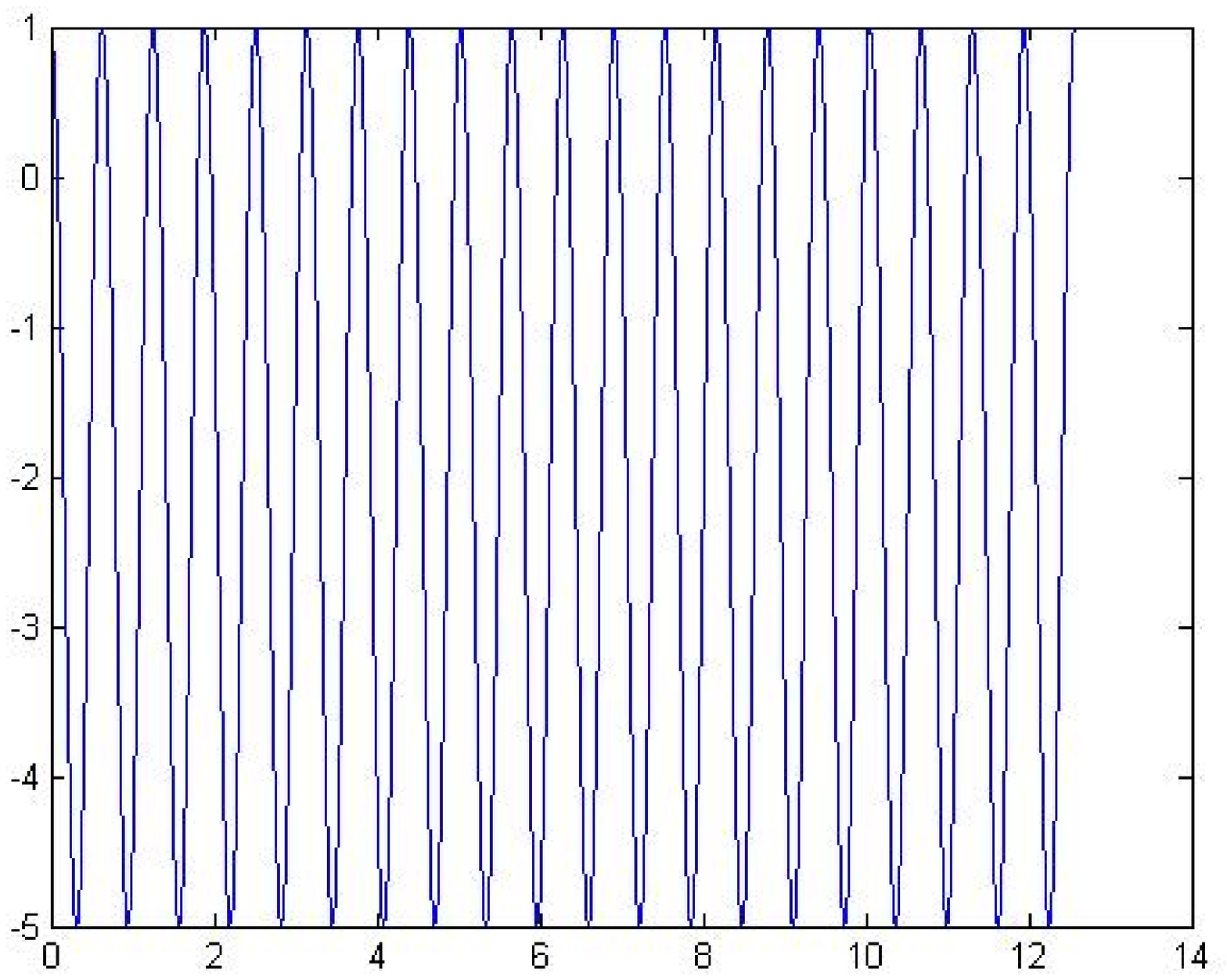}}
\subfigure[PerNoLin]{\includegraphics[width=0.32\textwidth]{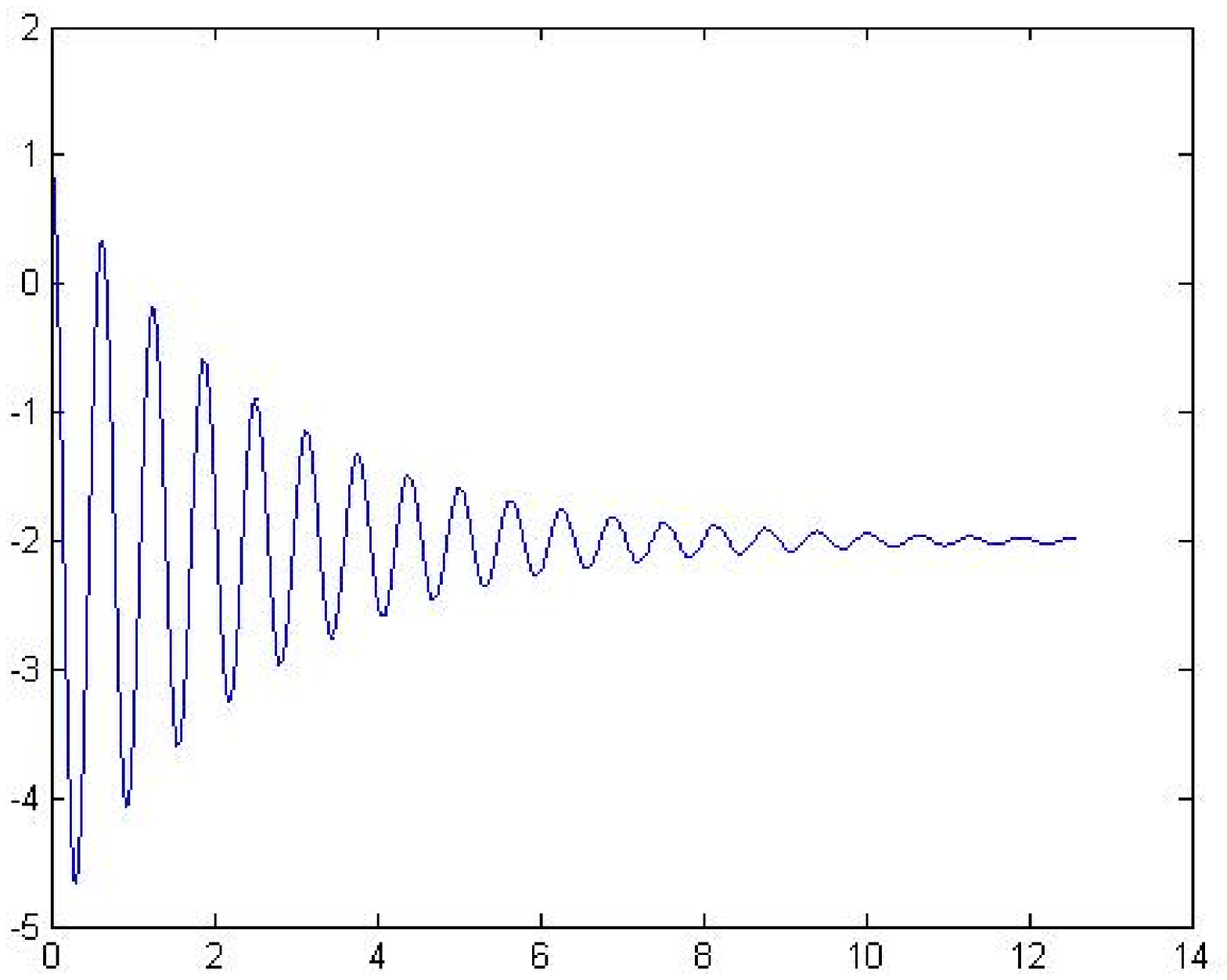}}
\subfigure[StiffLin]{\includegraphics[width=0.32\textwidth]{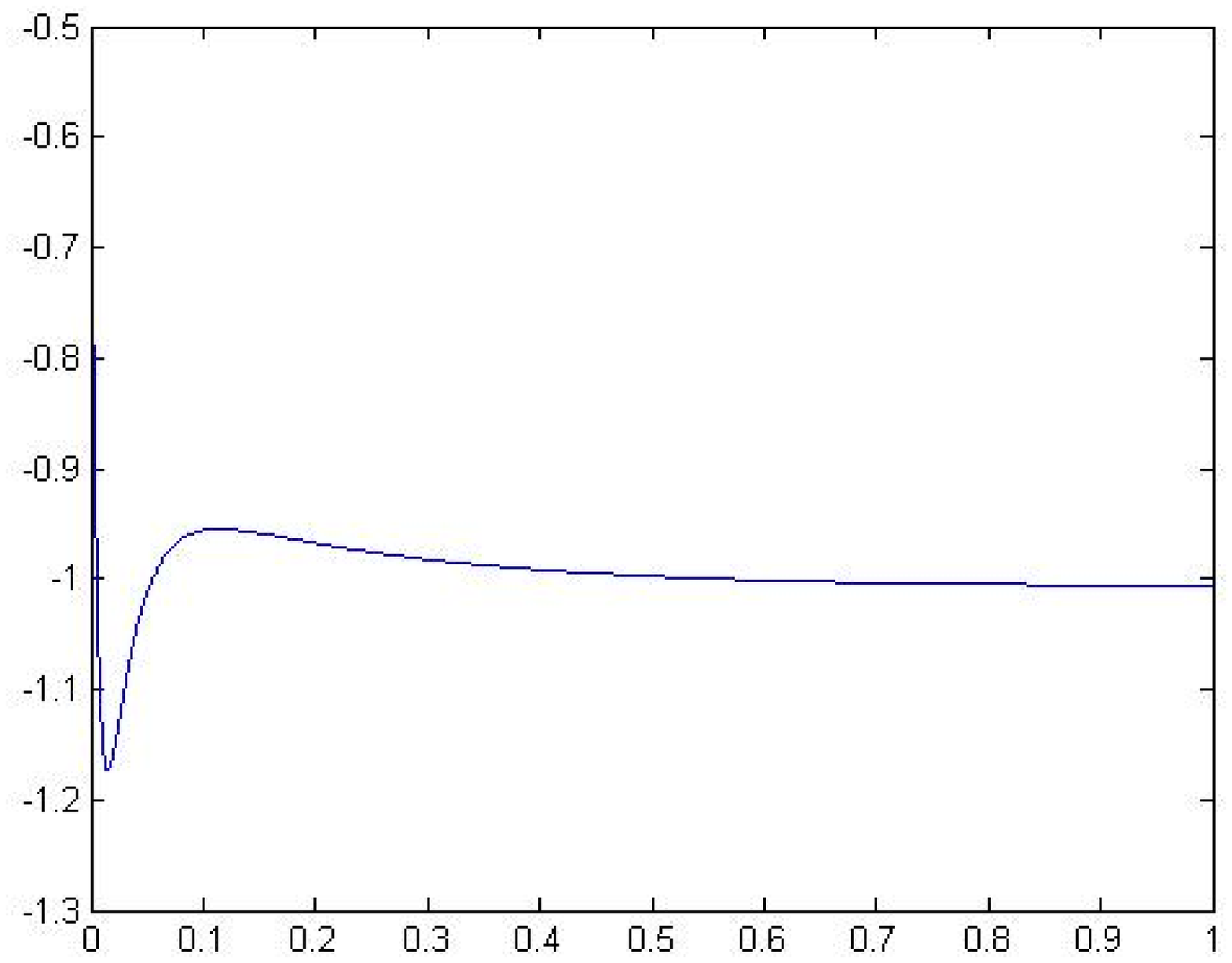}}\\
\subfigure[StiffNoLin]{\includegraphics[width=0.32\textwidth]{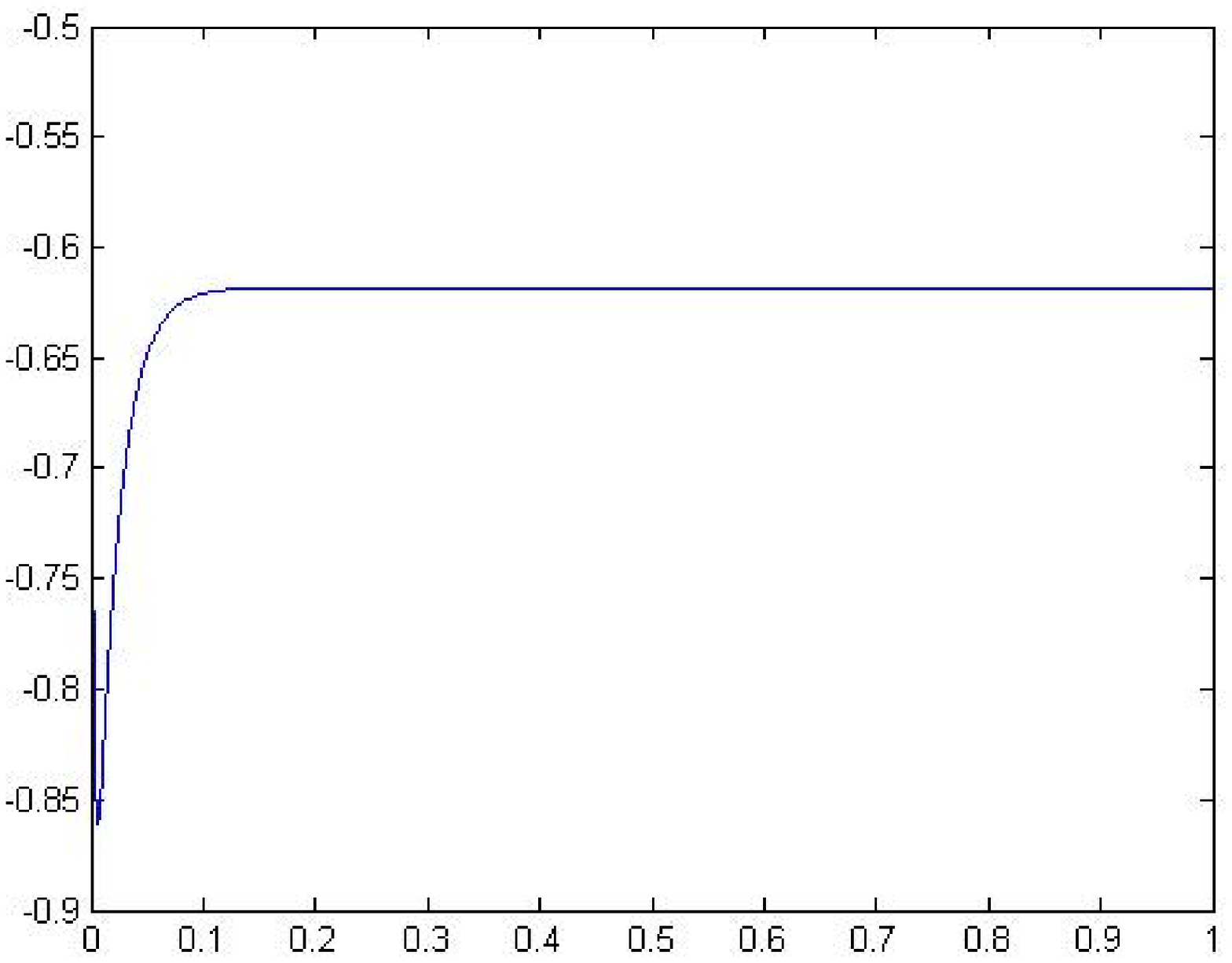}}
\subfigure[fpu]{\includegraphics[width=0.32\textwidth]{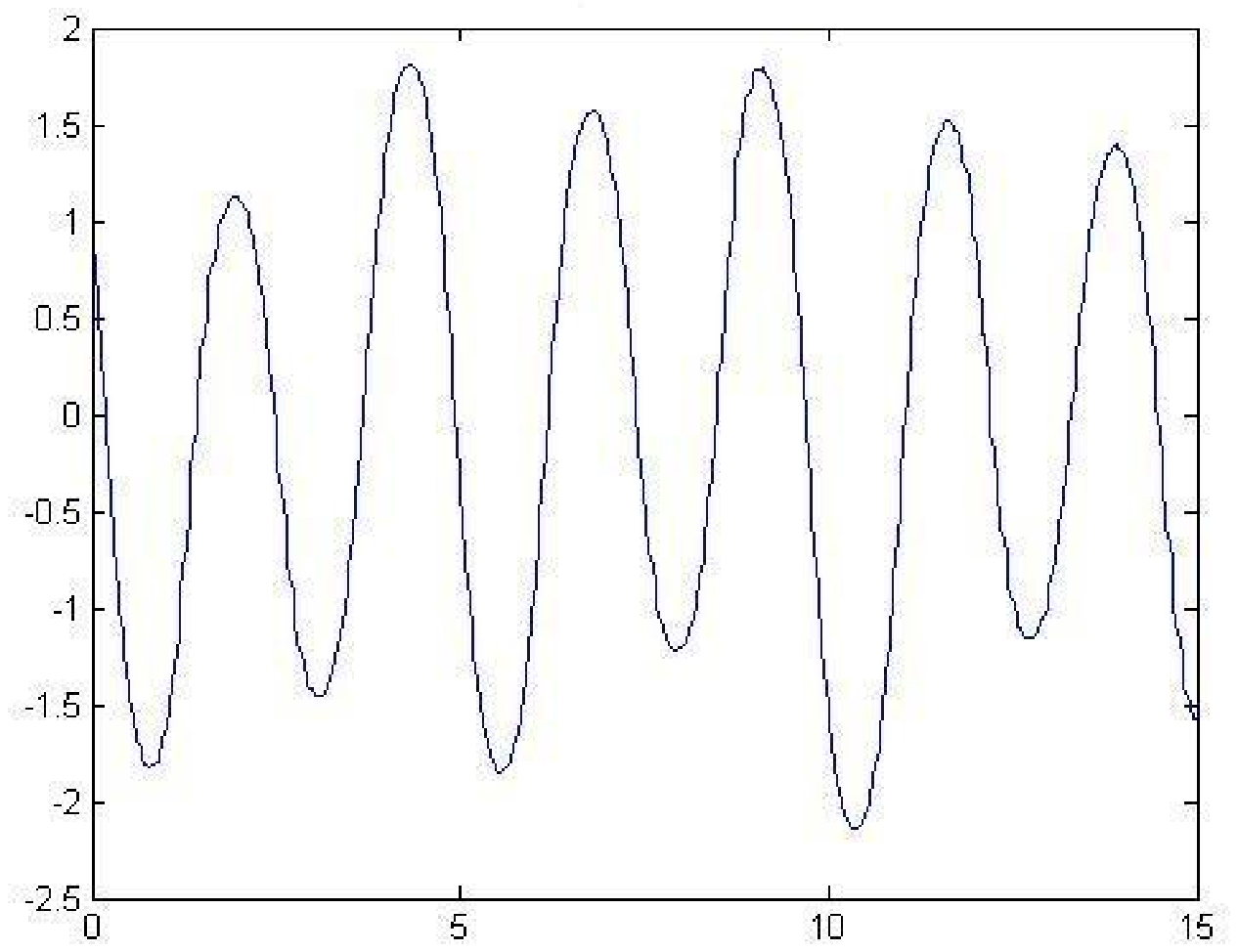}}
\subfigure[bruss]{\includegraphics[width=0.32\textwidth]{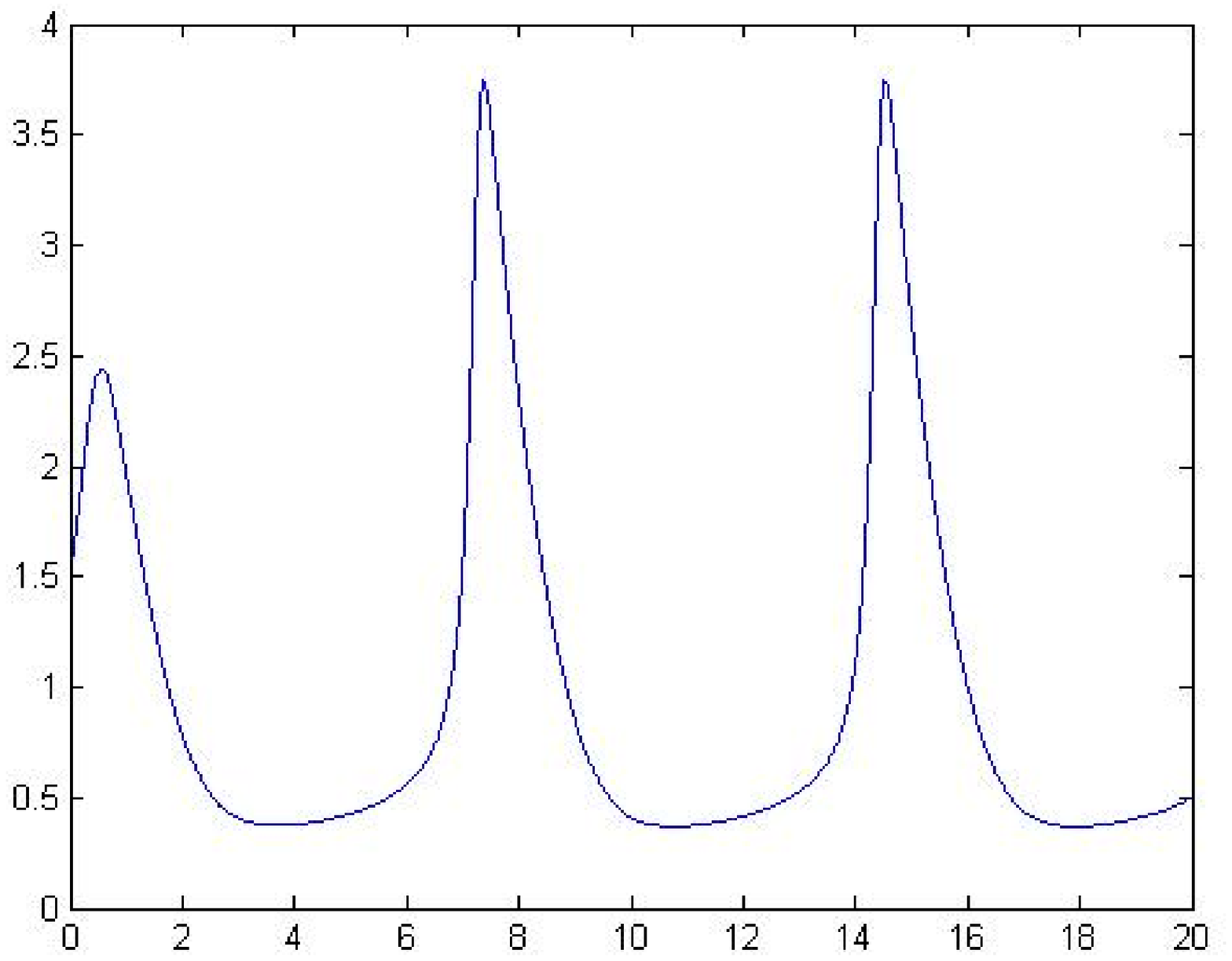}}\\
\subfigure[rigid]{\includegraphics[width=0.32\textwidth]{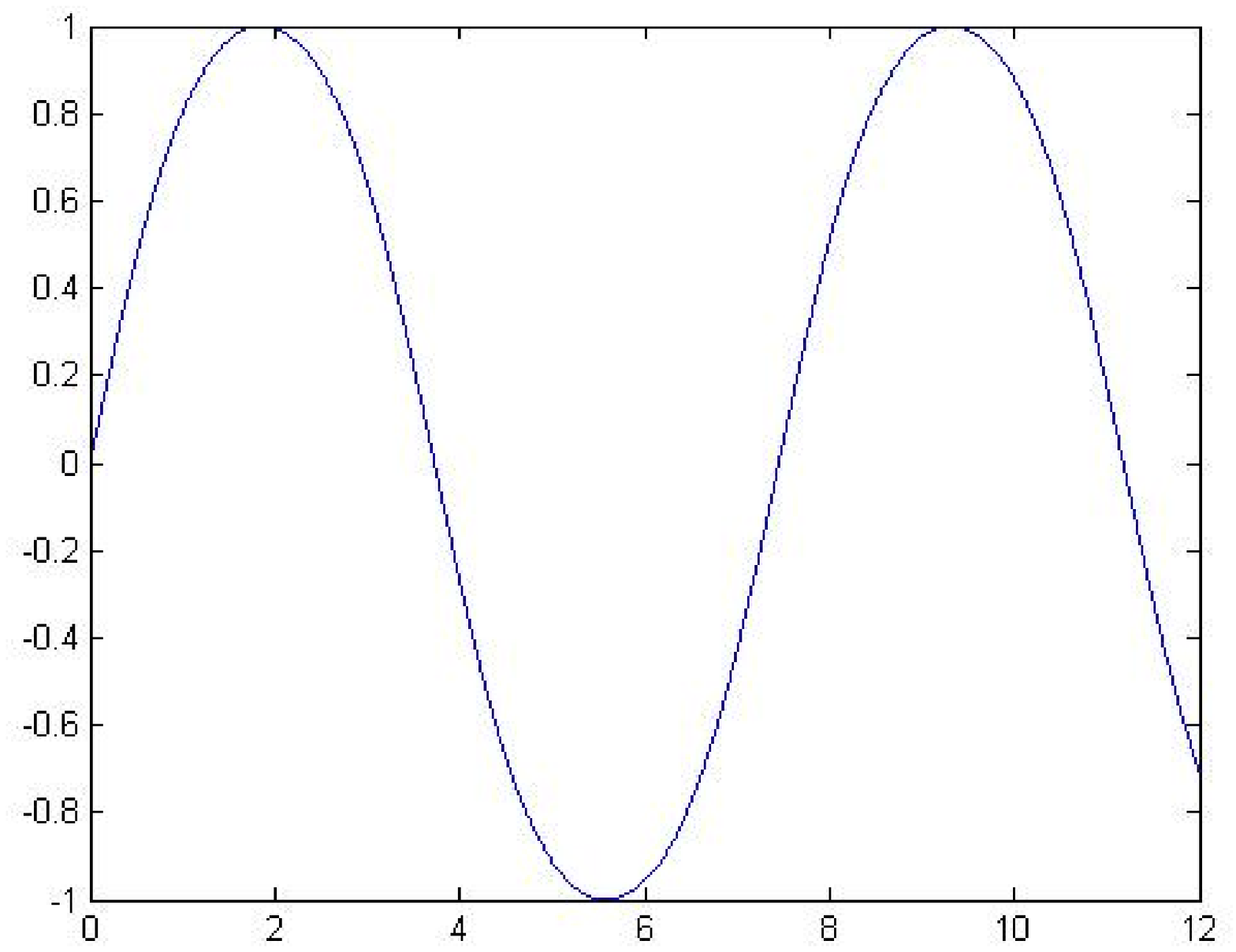}}
\subfigure[chm]{\includegraphics[width=0.32\textwidth]{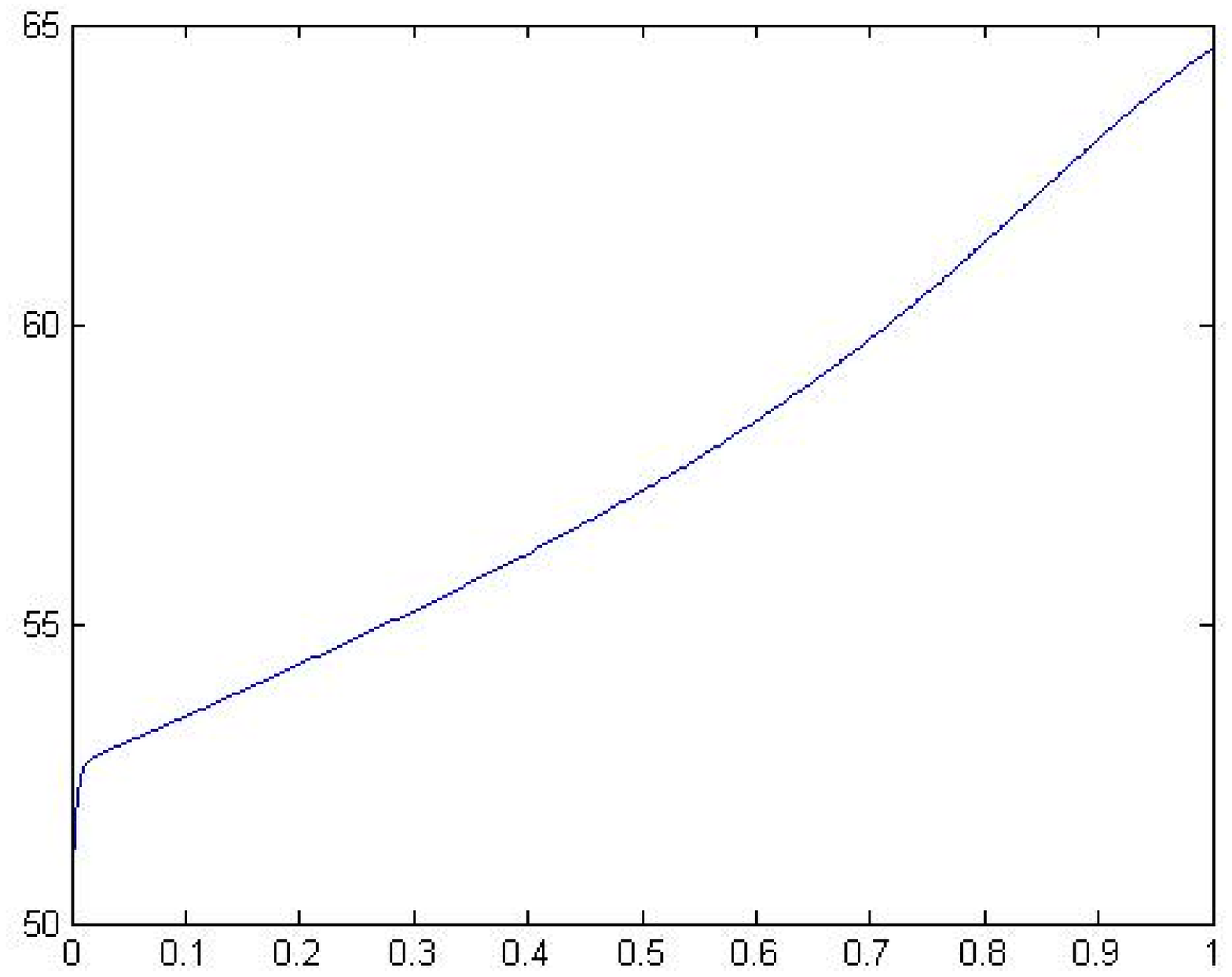}}
\subfigure[vdp1]{\includegraphics[width=0.32\textwidth]{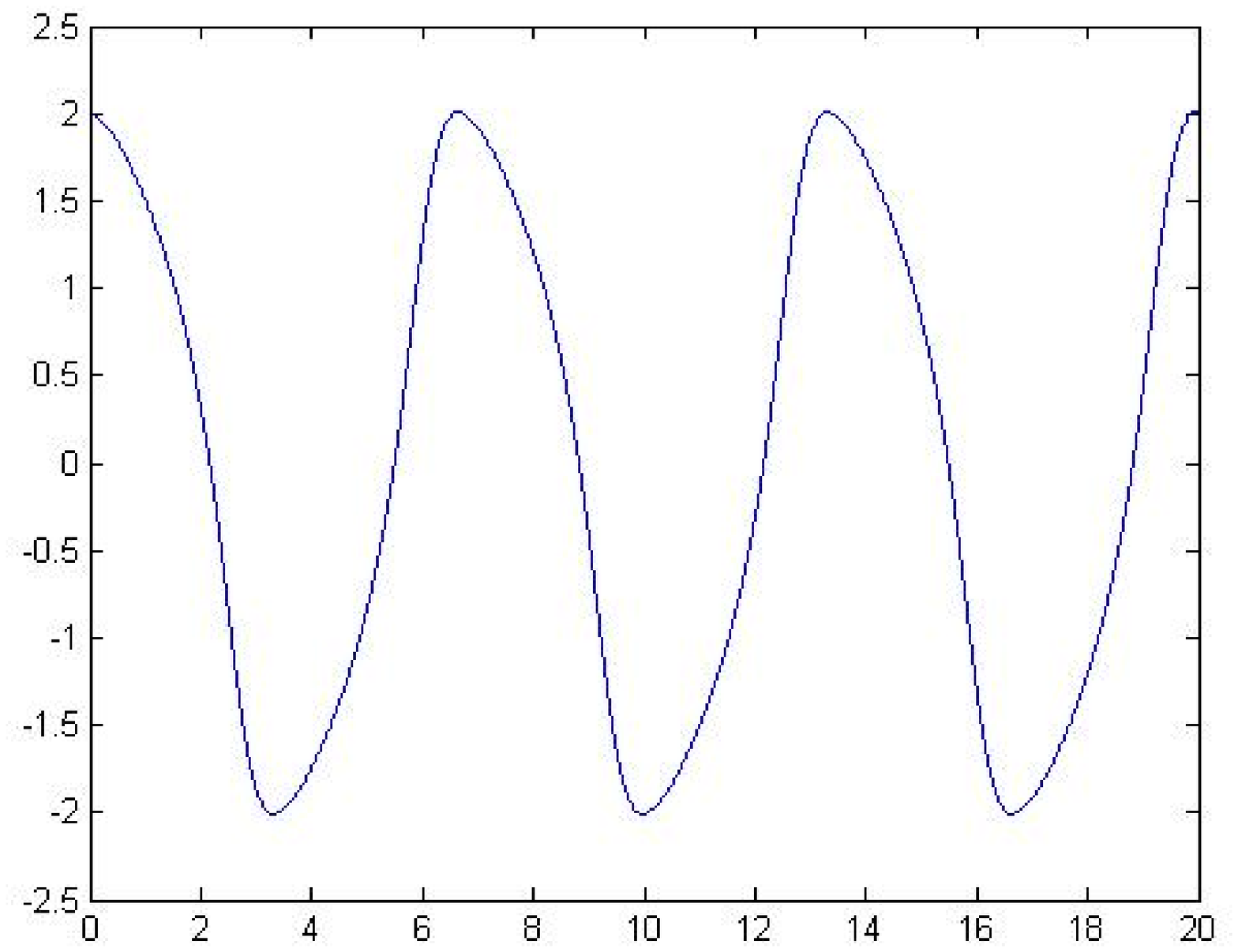}}\\
\subfigure{\includegraphics[width=0.32\textwidth]{blank.eps}}
\subfigure[vdp100]{\includegraphics[width=0.32\textwidth]{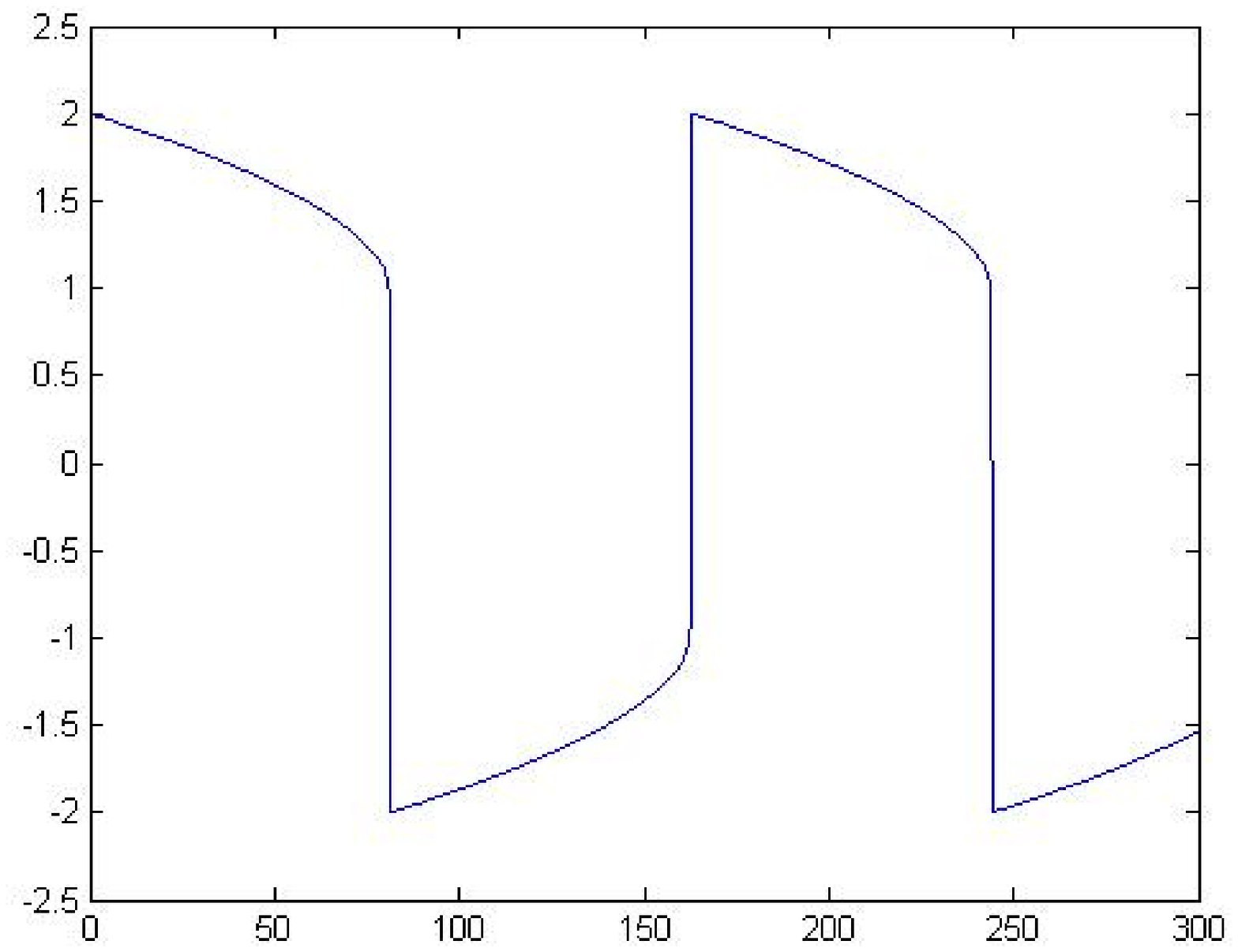}}
\caption{Path of the first components of the solution in each
example\label{test equations}.}
\end{figure}

As illustration, Figure \ref{test equations} shows the first
component of the solution of each example, which will be
consecutively named as \textit{PerLin, PerNoLin, StiffLin,
StiffNoLin, fpu, bruss, rigid, chm, vdp1 and vdp100}.

\subsection{Simulation A: integration over same time partition}

This simulation is designed to compare the accuracy of the order-$5$
formulas of the codes LLDP$45$ and ode$45$ over identical time
partitions. First, the ode$45$ code integrates all the examples with
the crude tolerances $RTol=10^{-3}$ and $ATol=10^{-6}$. This
defined, for each example, a time partition $ (t)_{h} $ over which
the order-$5$ formula of the LLDP$45$ code is evaluated as well.
That is, the formula
\begin{equation}
\mathbf{y}_{n+1}=\mathbf{y}_{n}+\mathbf{u}_{s}+h_{n}\sum_{j=1}^{s}b_{j}%
\mathbf{k}_{j}  \label{Order-5 LLRK formula}
\end{equation}
with $\mathbf{u}_{j}=\mathbf{L}(\mathbf{P}_{3,3}(2^{-\kappa
_{j}}\mathbf{D} _{n}c_{j}h_{n}))^{2^{\kappa _{j}}}\mathbf{r}$ for
the LLDP$45$ code. Tables \ref{Table I ODE} and \ref{Table II ODE}
present, respectively, the Relative Error (\ref{Relative Error}) of
the order-$5$ formula of each code in the integration of the four
semilinear and six nonlinear examples defined above. The number of
accepted time steps is shown as well. This comparison is repeated
twice but with the mild and refined tolerances
$RTol=10^{-6}$,$ATol=10^{-9}$ and $RTol=10^{-9}$,$ATol=10^{-12}$.
The results are also shown in Tables \ref{Table I ODE} and
\ref{Table II ODE}.
\newline

\begin{table}[h] \centering%
{\scriptsize $%
\begin{tabular}{|c|c|c|c|c|c}
\hline
Example & Tol & $%
\begin{array}{c}
\text{Time} \\
\text{steps}%
\end{array}%
$ & $%
\begin{array}{c}
\text{Relative Error} \\
\text{DP formulas}
\end{array}%
$ & $%
\begin{array}{c}
\text{Relative Error} \\
\text{LLDP formulas}
\end{array}$\\ \hline PerLin  &
$\begin{array}{c}
  \text{Crude} \\
  \text{Mild} \\
  \text{Refined}
\end{array}$
 & \multicolumn{1}{|c|}{$\begin{array}{c}
147 \\
598 \\
2394
\end{array}$} & $\begin{array}{c}
1.01 \\
2.6\times 10^{-2} \\
2.7\times10^{-5}
\end{array}$
  & $\begin{array}{c}
  2.0\times 10^{-6} \\
  3.2\times 10^{-7} \\
  1.3\times 10^{-6}
\end{array}$\\\cline{2-5} \hline PerNoLin &
$\begin{array}{c}
  \text{Crude} \\
  \text{Mild} \\
  \text{Refined}
\end{array}$ & \multicolumn{1}{|c|}{$\begin{array}{c}
105 \\
411 \\
1634
\end{array}$} & $\begin{array}{c}
  1.3\times10^{-2} \\
  2.6\times10^{-5} \\
  1.2\times10^{-7}
  \end{array}$  & $\begin{array}{c}
  4.9\times10^{-5} \\
  6.9\times 10^{-8} \\
  1.4\times10^{-9}
  \end{array}$ \\ \cline{2-5} \hline
 StiffLin  &
$\begin{array}{c}
  \text{Crude} \\
  \text{Mild} \\
  \text{Refined}
\end{array}$ & \multicolumn{1}{|c|}{$\begin{array}{c}
60 \\
78 \\
173
\end{array}$} &
$\begin{array}{c}
  1.1\times10^{-3} \\
  1.1\times10^{-6} \\
  8.0\times10^{-10}
  \end{array}$  &
$\begin{array}{c}
  2.7\times10^{-12} \\
  2.7\times 10^{-12} \\
  2.7\times10^{-12}
  \end{array}$ \\ \cline{2-5}  \hline
StiffNoLin & $\begin{array}{c}
  \text{Crude} \\
  \text{Mild} \\
  \text{Refined}
\end{array}$ &
\multicolumn{1}{|c|}{$\begin{array}{c}
104 \\
133 \\
294
\end{array}$} & $\begin{array}{c}
  1.4\times10^{-2} \\
  1.5\times 10^{-5} \\
  1.4\times10^{-8}
  \end{array}$ & $\begin{array}{c}
  9.7\times10^{-5} \\
  6.8\times 10^{-8} \\
  1.3\times10^{-8}
  \end{array}$
 \\ \cline{2-5}\hline
\end{tabular}%
$}%
\caption{Relative error of the order-5 formula of each code when
integrate
the semilinear examples over identical time partition.}\label{Table I ODE}%
\end{table}%

Note that, in general, the ode$45$ code is able to adequately
integrate the test equations with the tree specified tolerances.
Exceptions are the highly oscillatory \textit{fpu} equation and the
moderate stiff equation \textit{vpd100} at crude tolerances, for
which the relative error is high or unacceptable, respectively.
Observe that, the order-$5$ locally linearized formula (\ref{Order-5
LLRK formula}) is able to integrate the first equation with an
adequate relative error, but fail to integrate the second one on the
time partition generated by the ode$45$ code. In this last case, the
Padé algorithm fails to compute the exponential matrix at some point
of the mentioned time partition and, because of that, the place
corresponding to this information in Table \ref{Table II ODE} is
empty.
\newline
\newline

\begin{table}[h] \centering%
{\scriptsize $%
\begin{tabular}{|c|c|c|c|c|}
\hline Example & Tol & $%
\begin{array}{c}
\text{Time} \\
\text{steps}%
\end{array}%
$ &$
\begin{array}{c}
\text{Relative Error} \\
\text{DP formulas}
\end{array}%
$ &$ \begin{array}{c}
\text{Relative Error} \\
\text{LLDP formulas}
\end{array}$\\ \hline fpu  &
$\begin{array}{c}
  \text{Crude} \\
  \text{Mild} \\
  \text{Refined}
\end{array}$ & $\begin{array}{c}
  964 \\
  4474 \\
  19190
\end{array}$ & $\begin{array}{c}
  1.9\times 10^{2}\\
  8.1 \\
  1.7\times 10^{-2}
\end{array}$
 & $\begin{array}{c}
  1.5\times 10^{-2} \\
  2.9\times 10^{-3} \\
  1.7\times 10^{-2}
\end{array}$ \\
\hline rigid  & $\begin{array}{c}
  \text{Crude} \\
  \text{Mild} \\
  \text{Refined}
\end{array}$ &
$\begin{array}{c}
  19 \\
  66 \\
  256
\end{array}$ & $\begin{array}{c}
  2.7\times 10^{-2} \\
  7.5\times 10^{-5} \\
  2.0\times 10^{-7}
\end{array}$& $\begin{array}{c}
  1.5\times 10^{-3} \\
  4.0\times 10^{-6} \\
  1.8\times 10^{-8}
\end{array}$ \\ \hline chm &
$\begin{array}{c}
  \text{Crude} \\
  \text{Mild} \\
  \text{Refined}
\end{array}$ &
$\begin{array}{c}
  679 \\
  723 \\
  1520
\end{array}$& $\begin{array}{c}
  1.1\times 10^{-3} \\
  1.1\times 10^{-6} \\
  1.2\times 10^{-8}
\end{array}$   & $\begin{array}{c}
  5.5\times 10^{-7} \\
  2.5\times 10^{-7} \\
  1.2\times 10^{-8}
\end{array}$ \\ \hline bruss  & $\begin{array}{c}
  \text{Crude} \\
  \text{Mild} \\
  \text{Refined}
\end{array}$ &
$\begin{array}{c}
  46 \\
  148 \\
  558
\end{array}$ &
$\begin{array}{c}
  7.7\times 10^{-2} \\
  8.7\times 10^{-6} \\
  1.5\times 10^{-8}
\end{array}$  & $\begin{array}{c}
  2.4\times 10^{-2} \\
  3.5\times 10^{-7} \\
  1.2\times 10^{-9}
\end{array}$
 \\ \hline vdp1 &
$\begin{array}{c}
  \text{Crude} \\
  \text{Mild} \\
  \text{Refined}
\end{array}$ &
$\begin{array}{c}
  59 \\
  204 \\
  785
\end{array}$ & $\begin{array}{c}
  2.0 \\
  2.8\times 10^{-4} \\
  5.6\times 10^{-7}
\end{array}$ & $\begin{array}{c}
  0.14 \\
  1.5\times 10^{-5} \\
  3.1\times 10^{-8}
\end{array}$
 \\ \hline vdp100 &
$\begin{array}{c}
  \text{Crude} \\
  \text{Mild} \\
  \text{Refined}
\end{array}$ & $\begin{array}{c}
  16916 \\
  17516 \\
  31254
\end{array}$ & $\begin{array}{c}
  1.9\times 10^{4} \\
  0.42 \\
  1.2\times 10^{-3}
\end{array}$ & $\begin{array}{c}
   - \\
  9.2\times 10^{-2} \\
  8.1\times 10^{-4}
\end{array}$\\ \hline
\end{tabular}%
$}%
\caption{Relative error of the order-5 formula of each code when
integrate
the nonlinear examples over identical time partition.}\label{Table II ODE}%
\end{table}%

\subsection{Simulation B: integration with same tolerance}

This simulation is designed to compare the performance the codes
LLDP$45$ and ode$45$ with the same tolerances. As a difference with
Simulation A, here each codes use a different time partition defined
by their own adaptive strategy.

Tables \ref{Table III ODE} and \ref{Table IV ODE} summarize the
results of each code in the integration of each example for the
three sets of tolerances specified above. The column\ "Time" in
these tables presents the relative overall time of each numerical
scheme with respect to that of the ode45 code on the whole interval
$[t_{0},T]$. This overall time ratio works as simple indicator to
compare the total computational cost of each code. In addition, the
tables show the number of accepted and failed steps, the number of
evaluations of $\mathbf{f}$ and $ \mathbf{f}_{x}$, and the number of
exponential matrices computed in the integration of each example.
\newline
\newline
\newline
\newline

\begin{table}[H]
\centering{\scriptsize $%
\begin{tabular}{|c|c|c|c|c|c|c|c|c|}
\hline
Example & Code & Tol & $%
\begin{array}{c}
\text{Time} \\
\text{steps}%
\end{array}%
$ & $%
\begin{array}{c}
\text{Failed} \\
\text{steps}%
\end{array}%
$ & $f$\  & exp($f_{x})$ & Time & $%
\begin{array}{c}
\text{Relative} \\
\text{Error}%
\end{array}%
$ \\ \hline
PerLin & \multicolumn{1}{|l|}{ode45} & $%
\begin{array}{c}
\text{Crude} \\
\text{Mild} \\
\text{Refined}%
\end{array}%
$ & $%
\begin{array}{c}
147 \\
598 \\
2394%
\end{array}%
$ & $%
\begin{array}{c}
0 \\
0 \\
0%
\end{array}%
$ & $%
\begin{array}{c}
883 \\
3589 \\
14365%
\end{array}%
$ & $%
\begin{array}{c}
0 \\
0 \\
0%
\end{array}%
$ & $%
\begin{array}{c}
1 \\
1 \\
1%
\end{array}%
$ & $%
\begin{array}{c}
10.2 \\
2.6\times 10^{-2} \\
2.7\times 10^{-5}%
\end{array}%
$ \\ \cline{2-9}
& \multicolumn{1}{|l|}{LLDP$45$} & $%
\begin{array}{c}
\text{Crude} \\
\text{Mild} \\
\text{Refined}%
\end{array}%
$ & $%
\begin{array}{c}
14 \\
14 \\
15%
\end{array}%
$ & $%
\begin{array}{c}
0 \\
0 \\
0%
\end{array}%
$ & $%
\begin{array}{c}
85 \\
85 \\
91%
\end{array}%
$ & $%
\begin{array}{c}
14 \\
14 \\
15%
\end{array}%
$ & $%
\begin{array}{c}
0.27 \\
0.08 \\
0.01%
\end{array}%
$ & $%
\begin{array}{c}
2.0\times 10^{-9} \\
3.0\times 10^{-9} \\
2.0\times 10^{-9}%
\end{array}%
$ \\ \hline
PerNoLin & \multicolumn{1}{|l|}{ode45} & $%
\begin{array}{c}
\text{Crude} \\
\text{Mild} \\
\text{Refined}%
\end{array}%
$ & $%
\begin{array}{c}
105 \\
411 \\
1634%
\end{array}%
$ & $%
\begin{array}{c}
0 \\
0 \\
0%
\end{array}%
$ & $%
\begin{array}{c}
631 \\
2467 \\
9805%
\end{array}%
$ & $%
\begin{array}{c}
0 \\
0 \\
0%
\end{array}%
$ & $%
\begin{array}{c}
1 \\
1 \\
1%
\end{array}%
$ & $%
\begin{array}{c}
1.3\times 10^{-2} \\
2.6\times 10^{-5} \\
1.2\times 10^{-7}%
\end{array}%
$ \\ \cline{2-9}
& \multicolumn{1}{|l|}{LLDP$45$} & $%
\begin{array}{c}
\text{Crude} \\
\text{Mild} \\
\text{Refined}%
\end{array}%
$ & $%
\begin{array}{c}
42 \\
137 \\
534%
\end{array}%
$ & $%
\begin{array}{c}
0 \\
0 \\
0%
\end{array}%
$ & $%
\begin{array}{c}
253 \\
823 \\
3205%
\end{array}%
$ & $%
\begin{array}{c}
42 \\
137 \\
534%
\end{array}%
$ & $%
\begin{array}{c}
0.74 \\
0.62 \\
0.57%
\end{array}%
$ & $%
\begin{array}{c}
2.2\times 10^{-3} \\
3.6\times 10^{-6} \\
2.1\times 10^{-9}%
\end{array}%
$ \\ \hline
StiffLin & \multicolumn{1}{|l|}{ode45} & $%
\begin{array}{c}
\text{Crude} \\
\text{Mild} \\
\text{Refined}%
\end{array}%
$ & $%
\begin{array}{c}
60 \\
78 \\
172%
\end{array}%
$ & $%
\begin{array}{c}
6 \\
1 \\
6%
\end{array}%
$ & $%
\begin{array}{c}
397 \\
475 \\
1069%
\end{array}%
$ & $%
\begin{array}{c}
0 \\
0 \\
0%
\end{array}%
$ & $%
\begin{array}{c}
1 \\
1 \\
1%
\end{array}%
$ & $%
\begin{array}{c}
1.1\times 10^{-3} \\
1.1\times 10^{-6} \\
8.0\times 10^{-10}%
\end{array}%
$ \\ \cline{2-9}
& \multicolumn{1}{|l|}{LLDP$45$} & $%
\begin{array}{c}
\text{Crude} \\
\text{Mild} \\
\text{Refined}%
\end{array}%
$ & $%
\begin{array}{c}
14 \\
14 \\
15%
\end{array}%
$ & $%
\begin{array}{c}
0 \\
0 \\
0%
\end{array}%
$ & $%
\begin{array}{c}
85 \\
85 \\
91%
\end{array}%
$ & $%
\begin{array}{c}
14 \\
14 \\
15%
\end{array}%
$ & $%
\begin{array}{c}
0.33 \\
0.34 \\
0.15%
\end{array}%
$ & $%
\begin{array}{c}
2.5\times 10^{-12} \\
2.3\times 10^{-12} \\
2.3\times 10^{-12}%
\end{array}%
$ \\ \hline
StiffNoLin & \multicolumn{1}{|l|}{ode45} & $%
\begin{array}{c}
\text{Crude} \\
\text{Mild} \\
\text{Refined}%
\end{array}%
$ & $%
\begin{array}{c}
104 \\
133 \\
294%
\end{array}%
$ & $%
\begin{array}{c}
4 \\
5 \\
2%
\end{array}%
$ & $%
\begin{array}{c}
649 \\
829 \\
1777%
\end{array}%
$ & $%
\begin{array}{c}
0 \\
0 \\
0%
\end{array}%
$ & $%
\begin{array}{c}
1 \\
1 \\
1%
\end{array}%
$ & $%
\begin{array}{c}
1.4\times 10^{-2} \\
1.5\times 10^{-5} \\
1.4\times 10^{-8}%
\end{array}%
$ \\ \cline{2-9}
& \multicolumn{1}{|l|}{LLDP$45$} & $%
\begin{array}{c}
\text{Crude} \\
\text{Mild} \\
\text{Refined}%
\end{array}%
$ & $%
\begin{array}{c}
21 \\
43 \\
132%
\end{array}%
$ & $%
\begin{array}{c}
0 \\
0 \\
2%
\end{array}%
$ & $%
\begin{array}{c}
127 \\
259 \\
805%
\end{array}%
$ & $%
\begin{array}{c}
21 \\
43 \\
134%
\end{array}%
$ & $%
\begin{array}{c}
0.32 \\
0.53 \\
0.68%
\end{array}%
$ & $%
\begin{array}{c}
8.0\times 10^{-4} \\
1.6\times 10^{-6} \\
9.2\times 10^{-9}%
\end{array}%
$ \\ \hline
\end{tabular}%
$} \caption{Code performance in the integration of the semilinear
examples with the same tolerances.} \label{Table III ODE}
\end{table}

\subsection{Simulation C: integration with similar accuracy}

In this type of simulation, the tolerances $RTol$ and $ATol$ of the
LLDP$45$ code is changed until its relative error in the integration
of each example achieves similar value to that corresponding to the
code ode$45$. This simulation is carried out three times, changing
the tolerances of the ode$45$ from the crude values to the refined
values specified above. Tables \ref {Table V ODE} and \ref{Table VI
ODE} summarize the performance of each code in the integration of
each example. As in the previous two tables, this includes the
relative overall time, number of accepted and failed steps, the
number of evaluations of $\mathbf{f}$ and $ \mathbf{f}_{x}$, and the
number of exponential matrices computed in the integration of each
example.
\newline
\newline
\newline
\newline
\newline
\newline
\newline
\newline

\begin{table}[H]
\centering{\scriptsize $%
\begin{tabular}{|c|c|c|c|c|c|c|c|c|}
\hline
Example & Code & Tol & $%
\begin{array}{c}
\text{Time} \\
\text{steps}%
\end{array}%
$ & $%
\begin{array}{c}
\text{Failed} \\
\text{steps}%
\end{array}%
$ & $f$\  & exp($f_{x})$ & Time & $%
\begin{array}{c}
\text{Relative} \\
\text{Error}%
\end{array}%
$ \\ \hline
\multicolumn{1}{|c|}{fpu} & \multicolumn{1}{|l|}{ode45} & $%
\begin{array}{c}
\text{Crude} \\
\text{Mild} \\
\text{Refined}%
\end{array}%
$ & $%
\begin{array}{c}
964 \\
4474 \\
19190%
\end{array}%
$ & $%
\begin{array}{c}
2 \\
60 \\
45%
\end{array}%
$ & $%
\begin{array}{c}
5797 \\
27205 \\
115411%
\end{array}%
$ & $%
\begin{array}{c}
0 \\
0 \\
0%
\end{array}%
$ & $%
\begin{array}{c}
1 \\
1 \\
1%
\end{array}%
$ & $%
\begin{array}{c}
1.9\times 10^{2} \\
8.1 \\
1.7\times 10^{-2}%
\end{array}%
$ \\ \cline{2-9}
& \multicolumn{1}{|l|}{LLDP$45$} & $%
\begin{array}{c}
\text{Crude} \\
\text{Mild} \\
\text{Refined}%
\end{array}%
$ & $%
\begin{array}{c}
377 \\
1496 \\
6021%
\end{array}%
$ & $%
\begin{array}{c}
48 \\
125 \\
86%
\end{array}%
$ & $%
\begin{array}{c}
2551 \\
9727 \\
36643%
\end{array}%
$ & $%
\begin{array}{c}
425 \\
1621 \\
6107%
\end{array}%
$ & $%
\begin{array}{c}
0.81 \\
0.67 \\
0.49%
\end{array}%
$ & $%
\begin{array}{c}
17.4\\
2.0\times 10^{-2} \\
1.7\times 10^{-2}%
\end{array}%
$ \\ \hline
rigid & \multicolumn{1}{|l|}{ode45} & $%
\begin{array}{c}
\text{Crude} \\
\text{Mild} \\
\text{Refined}%
\end{array}%
$ & $%
\begin{array}{c}
19 \\
66 \\
256%
\end{array}%
$ & $%
\begin{array}{c}
2 \\
4 \\
1%
\end{array}%
$ & $%
\begin{array}{c}
127 \\
421 \\
1543%
\end{array}%
$ & $%
\begin{array}{c}
0 \\
0 \\
0%
\end{array}%
$ & $%
\begin{array}{c}
1 \\
1 \\
1%
\end{array}%
$ & $%
\begin{array}{c}
2.7\times 10^{-2} \\
7.5\times 10^{-5} \\
2.0\times 10^{-7}%
\end{array}%
$ \\ \cline{2-9}
& \multicolumn{1}{|l|}{LLDP$45$} & $%
\begin{array}{c}
\text{Crude} \\
\text{Mild} \\
\text{Refined}%
\end{array}%
$ & $%
\begin{array}{c}
16 \\
53 \\
201%
\end{array}%
$ & $%
\begin{array}{c}
0 \\
5 \\
0%
\end{array}%
$ & $%
\begin{array}{c}
97 \\
349 \\
1207%
\end{array}%
$ & $%
\begin{array}{c}
16 \\
58 \\
201%
\end{array}%
$ & $%
\begin{array}{c}
1.18 \\
1.55 \\
1.48%
\end{array}%
$ & $%
\begin{array}{c}
3.3\times 10^{-3} \\
8.6\times 10^{-6} \\
3.1\times 10^{-8}%
\end{array}%
$ \\ \hline
chm & \multicolumn{1}{|l|}{ode45} & $%
\begin{array}{c}
\text{Crude} \\
\text{Mild} \\
\text{Refined}%
\end{array}%
$ & $%
\begin{array}{c}
679 \\
723 \\
1521%
\end{array}%
$ & $%
\begin{array}{c}
47 \\
16 \\
1%
\end{array}%
$ & $%
\begin{array}{c}
4357 \\
4435 \\
9133%
\end{array}%
$ & $%
\begin{array}{c}
0 \\
0 \\
0%
\end{array}%
$ & $%
\begin{array}{c}
1 \\
1 \\
1%
\end{array}%
$ & $%
\begin{array}{c}
1.1\times 10^{-3} \\
1.1\times 10^{-6} \\
1.2\times 10^{-8}%
\end{array}%
$ \\ \cline{2-9}
& \multicolumn{1}{|l|}{LLDP$45$} & $%
\begin{array}{c}
\text{Crude} \\
\text{Mild} \\
\text{Refined}%
\end{array}%
$ & $%
\begin{array}{c}
152 \\
357 \\
859%
\end{array}%
$ & $%
\begin{array}{c}
1 \\
2 \\
57%
\end{array}%
$ & $%
\begin{array}{c}
919 \\
2155 \\
5497%
\end{array}%
$ & $%
\begin{array}{c}
153 \\
359 \\
916%
\end{array}%
$ & $%
\begin{array}{c}
0.43 \\
1.05 \\
1.18%
\end{array}%
$ & $%
\begin{array}{c}
8.4\times 10^{-4} \\
9.2\times 10^{-7} \\
1.2\times 10^{-8}%
\end{array}%
$ \\ \hline
bruss & \multicolumn{1}{|l|}{ode45} & $%
\begin{array}{c}
\text{Crude} \\
\text{Mild} \\
\text{Refined}%
\end{array}%
$ & $%
\begin{array}{c}
46 \\
148 \\
558%
\end{array}%
$ & $%
\begin{array}{c}
12 \\
13 \\
4%
\end{array}%
$ & $%
\begin{array}{c}
349 \\
967 \\
3373%
\end{array}%
$ & $%
\begin{array}{c}
0 \\
0 \\
0%
\end{array}%
$ & $%
\begin{array}{c}
1 \\
1 \\
1%
\end{array}%
$ & $%
\begin{array}{c}
7.7\times 10^{-2} \\
8.7\times 10^{-6} \\
1.5\times 10^{-8}%
\end{array}%
$ \\ \cline{2-9}
& \multicolumn{1}{|l|}{LLDP$45$} & $%
\begin{array}{c}
\text{Crude} \\
\text{Mild} \\
\text{Refined}%
\end{array}%
$ & $%
\begin{array}{c}
36 \\
105 \\
396%
\end{array}%
$ & $%
\begin{array}{c}
7 \\
14 \\
11%
\end{array}%
$ & $%
\begin{array}{c}
259 \\
715 \\
2443%
\end{array}%
$ & $%
\begin{array}{c}
43 \\
119 \\
407%
\end{array}%
$ & $%
\begin{array}{c}
1.32 \\
1.47 \\
1.38%
\end{array}%
$ & $%
\begin{array}{c}
6.2\times 10^{-3} \\
5.4\times 10^{-6} \\
4.8\times 10^{-9}%
\end{array}%
$ \\ \hline
vdp1 & \multicolumn{1}{|l|}{ode45} & $%
\begin{array}{c}
\text{Crude} \\
\text{Mild} \\
\text{Refined}%
\end{array}%
$ & $%
\begin{array}{c}
59 \\
204 \\
785%
\end{array}%
$ & $%
\begin{array}{c}
10 \\
32 \\
19%
\end{array}%
$ & $%
\begin{array}{c}
415 \\
1417 \\
4825%
\end{array}%
$ & $%
\begin{array}{c}
0 \\
0 \\
0%
\end{array}%
$ & $%
\begin{array}{c}
1 \\
1 \\
1%
\end{array}%
$ & $%
\begin{array}{c}
2.24 \\
2.8\times 10^{-4} \\
5.7\times 10^{-7}%
\end{array}%
$ \\ \cline{2-9}
& \multicolumn{1}{|l|}{LLDP$45$} & $%
\begin{array}{c}
\text{Crude} \\
\text{Mild} \\
\text{Refined}%
\end{array}%
$ & $%
\begin{array}{c}
44 \\
162 \\
609%
\end{array}%
$ & $%
\begin{array}{c}
10 \\
38 \\
12%
\end{array}%
$ & $%
\begin{array}{c}
325 \\
1201 \\
3727%
\end{array}%
$ & $%
\begin{array}{c}
54 \\
200 \\
621%
\end{array}%
$ & $%
\begin{array}{c}
1.23 \\
1.45 \\
1.24%
\end{array}%
$ & $%
\begin{array}{c}
1.95 \\
5.8\times 10^{-5} \\
1.4\times 10^{-7}%
\end{array}%
$ \\ \hline
vdp100 & \multicolumn{1}{|l|}{ode45} & $%
\begin{array}{c}
\text{Crude} \\
\text{Mild} \\
\text{Refined}%
\end{array}%
$ & $%
\begin{array}{c}
16916 \\
17516 \\
31253%
\end{array}%
$ & $%
\begin{array}{c}
1074 \\
1540 \\
9%
\end{array}%
$ & $%
\begin{array}{c}
107941 \\
114337 \\
187573%
\end{array}%
$ & $%
\begin{array}{c}
0 \\
0 \\
0%
\end{array}%
$ & $%
\begin{array}{c}
1 \\
1 \\
1%
\end{array}%
$ & $%
\begin{array}{c}
1.9\times 10^{4} \\
0.41 \\
1.2\times 10^{-3}%
\end{array}%
$ \\ \cline{2-9}
& \multicolumn{1}{|l|}{LLDP$45$} & $%
\begin{array}{c}
\text{Crude} \\
\text{Mild} \\
\text{Refined}%
\end{array}%
$ & $%
\begin{array}{c}
3866 \\
7893 \\
19887%
\end{array}%
$ & $%
\begin{array}{c}
120 \\
19 \\
568%
\end{array}%
$ & $%
\begin{array}{c}
23917 \\
47473 \\
122731%
\end{array}%
$ & $%
\begin{array}{c}
3986 \\
7912 \\
20455%
\end{array}%
$ & $%
\begin{array}{c}
0.35 \\
0.69 \\
1.02%
\end{array}%
$ & $%
\begin{array}{c}
16.1 \\
2.1\times 10^{-3} \\
5.6\times 10^{-4}%
\end{array}%
$ \\ \hline
\end{tabular}%
$} \caption{Code performance in the integration of the nonlinear
examples with the same tolerances.} \label{Table IV ODE}
\end{table}

\begin{table}[H]
\centering{\scriptsize $%
\begin{tabular}{|c|c|c|c|c|c|c|c|c|}
\hline
Example & Code & Tol & $%
\begin{array}{c}
\text{Time} \\
\text{steps}%
\end{array}%
$ & $%
\begin{array}{c}
\text{Failed} \\
\text{steps}%
\end{array}%
$ & $f$\  & exp($f_{x})$ & Time & $%
\begin{array}{c}
\text{Relative} \\
\text{Error}%
\end{array}%
$ \\ \hline
PerLin & \multicolumn{1}{|l|}{ode45} & $%
\begin{array}{c}
\text{Crude} \\
\text{Mild} \\
\text{Refined}%
\end{array}%
$ & $%
\begin{array}{c}
147 \\
598 \\
2394%
\end{array}%
$ & $%
\begin{array}{c}
0 \\
0 \\
0%
\end{array}%
$ & $%
\begin{array}{c}
883 \\
3589 \\
14365%
\end{array}%
$ & $%
\begin{array}{c}
0 \\
0 \\
0%
\end{array}%
$ & $%
\begin{array}{c}
1 \\
1 \\
1%
\end{array}%
$ & $%
\begin{array}{c}
10.1 \\
2.6\times 10^{-2} \\
2.7\times 10^{-5}%
\end{array}%
$ \\ \cline{2-9}
& \multicolumn{1}{|l|}{LLDP$45$} & $%
\begin{array}{c}
100\times \text{Crude} \\
100\times \text{Crude} \\
100\times \text{Crude}%
\end{array}%
$ & $%
\begin{array}{c}
13 \\
13 \\
13%
\end{array}%
$ & $%
\begin{array}{c}
0 \\
0 \\
0%
\end{array}%
$ & $%
\begin{array}{c}
79 \\
79 \\
79%
\end{array}%
$ & $%
\begin{array}{c}
13 \\
13 \\
13%
\end{array}%
$ & $%
\begin{array}{c}
0.23 \\
0.06 \\
0.01%
\end{array}%
$ & $%
\begin{array}{c}
2.0\times 10^{-9} \\
2.0\times 10^{-9} \\
2.0\times 10^{-9}%
\end{array}%
$ \\ \hline
PerNoLin & \multicolumn{1}{|l|}{ode45} & $%
\begin{array}{c}
\text{Crude} \\
\text{Mild} \\
\text{Refined}%
\end{array}%
$ & $%
\begin{array}{c}
105 \\
411 \\
1634%
\end{array}%
$ & $%
\begin{array}{c}
0 \\
0 \\
0%
\end{array}%
$ & $%
\begin{array}{c}
361 \\
2467 \\
9805%
\end{array}%
$ & $%
\begin{array}{c}
0 \\
0 \\
0%
\end{array}%
$ & $%
\begin{array}{c}
1 \\
1 \\
1%
\end{array}%
$ & $%
\begin{array}{c}
1.2\times 10^{-2} \\
2.6\times 10^{-5} \\
1.2\times 10^{-7}%
\end{array}%
$ \\ \cline{2-9}
& \multicolumn{1}{|l|}{LLDP$45$} & $%
\begin{array}{c}
7.5\times \text{Crude} \\
7\times \text{Mild} \\
80\times \text{Refined}%
\end{array}%
$ & $%
\begin{array}{c}
32 \\
95 \\
224%
\end{array}%
$ & $%
\begin{array}{c}
0 \\
0 \\
0%
\end{array}%
$ & $%
\begin{array}{c}
193 \\
571 \\
1345%
\end{array}%
$ & $%
\begin{array}{c}
32 \\
95 \\
224%
\end{array}%
$ & $%
\begin{array}{c}
0.55 \\
0.44 \\
0.24%
\end{array}%
$ & $%
\begin{array}{c}
1.2\times 10^{-2} \\
1.5\times 10^{-5} \\
1.2\times 10^{-7}%
\end{array}%
$ \\ \hline
StiffLin & \multicolumn{1}{|l|}{ode45} & $%
\begin{array}{c}
\text{Crude} \\
\text{Mild} \\
\text{Refined}%
\end{array}%
$ & $%
\begin{array}{c}
60 \\
78 \\
172%
\end{array}%
$ & $%
\begin{array}{c}
6 \\
1 \\
6%
\end{array}%
$ & $%
\begin{array}{c}
397 \\
475 \\
1069%
\end{array}%
$ & $%
\begin{array}{c}
0 \\
0 \\
0%
\end{array}%
$ & $%
\begin{array}{c}
1 \\
1 \\
1%
\end{array}%
$ & $%
\begin{array}{c}
1.1\times 10^{-3} \\
1.1\times 10^{-6} \\
8.0\times 10^{-10}%
\end{array}%
$ \\ \cline{2-9}
& \multicolumn{1}{|l|}{LLDP$45$} & $%
\begin{array}{c}
100\times \text{Crude} \\
100\times \text{Crude} \\
100\times \text{Crude}%
\end{array}%
$ & $%
\begin{array}{c}
13 \\
13 \\
13%
\end{array}%
$ & $%
\begin{array}{c}
0 \\
0 \\
0%
\end{array}%
$ & $%
\begin{array}{c}
79 \\
79 \\
79%
\end{array}%
$ & $%
\begin{array}{c}
13 \\
13 \\
13%
\end{array}%
$ & $%
\begin{array}{c}
0.35 \\
0.28 \\
0.12%
\end{array}%
$ & $%
\begin{array}{c}
2.6\times 10^{-12} \\
2.6\times 10^{-12} \\
2.6\times 10^{-12}%
\end{array}%
$ \\ \hline
StiffNoLin & \multicolumn{1}{|l|}{ode45} & $%
\begin{array}{c}
\text{Crude} \\
\text{Mild} \\
\text{Refined}%
\end{array}%
$ & $%
\begin{array}{c}
104 \\
133 \\
294%
\end{array}%
$ & $%
\begin{array}{c}
4 \\
5 \\
2%
\end{array}%
$ & $%
\begin{array}{c}
649 \\
829 \\
1777%
\end{array}%
$ & $%
\begin{array}{c}
0 \\
0 \\
0%
\end{array}%
$ & $%
\begin{array}{c}
1 \\
1 \\
1%
\end{array}%
$ & $%
\begin{array}{c}
1.4\times 10^{-2} \\
1.5\times 10^{-5} \\
1.4\times 10^{-8}%
\end{array}%
$ \\ \cline{2-9}
& \multicolumn{1}{|l|}{LLDP$45$} & $%
\begin{array}{c}
9\times \text{Crude} \\
40\times \text{Mild} \\
8.45\times \text{Refined}%
\end{array}%
$ & $%
\begin{array}{c}
24 \\
28 \\
90%
\end{array}%
$ & $%
\begin{array}{c}
2 \\
0 \\
0%
\end{array}%
$ & $%
\begin{array}{c}
157 \\
169 \\
541%
\end{array}%
$ & $%
\begin{array}{c}
26 \\
28 \\
90%
\end{array}%
$ & $%
\begin{array}{c}
0.38 \\
0.35 \\
0.46%
\end{array}%
$ & $%
\begin{array}{c}
1.3\times 10^{-2} \\
1.1\times 10^{-5} \\
1.0\times 10^{-8}%
\end{array}%
$ \\ \hline
\end{tabular}%
$} \caption{Code performance in the integration of the semilinear
examples with similar accuracy.} \label{Table V ODE}
\end{table}

\subsection{Simulation D: evaluation of the dense output}

This simulation is designed to compare the accuracy of the
continuous formulas of the codes LLDP$45$ and ode$45$ over their
dense output. For this, both codes are applied first to each example
with the same crude tolerances $ RTol=10^{-3}$ and $ATol=10^{-6}$
but, the relative error of each code is now computed on its
respective dense output instead on the time partition $ (t)_{h}$
defined by the adaptive strategy. Tables \ref{Table VII ODE} and
\ref{Table VIII ODE} present these relative errors. The number of
accepted time steps and dense output times are also shown. This
comparison is repeated twice but with the mild and refined
tolerances $RTol=10^{-6}$,$ATol=10^{-9}$ and
$RTol=10^{-9}$,$ATol=10^{-12}$. The results are also shown in Tables
\ref{Table VII ODE} and \ref{Table VIII ODE}.
\newline
\newline
\newline
\newline
\newline
\newline
\newline
\newline

\begin{table}[H]
\centering{\scriptsize $%
\begin{tabular}{|c|c|c|c|c|c|c|c|c|}
\hline
Example & Code & Tol & $%
\begin{array}{c}
\text{Time} \\
\text{steps}%
\end{array}%
$ & $%
\begin{array}{c}
\text{Failed} \\
\text{steps}%
\end{array}%
$ & $f$\  & exp($f_{x})$ & Time & $%
\begin{array}{c}
\text{Relative} \\
\text{Error}%
\end{array}%
$ \\ \hline
fpu & \multicolumn{1}{|l|}{ode45} & $%
\begin{array}{c}
\text{Crude} \\
\text{Mild} \\
\text{Refined}%
\end{array}%
$ & $%
\begin{array}{c}
964 \\
4474 \\
19190%
\end{array}%
$ & $%
\begin{array}{c}
2 \\
60 \\
45%
\end{array}%
$ & $%
\begin{array}{c}
5797 \\
27205 \\
115411%
\end{array}%
$ & $%
\begin{array}{c}
0 \\
0 \\
0%
\end{array}%
$ & $%
\begin{array}{c}
1 \\
1 \\
1%
\end{array}%
$ & $%
\begin{array}{c}
1.9\times 10^{2} \\
8.1 \\
1.7\times 10^{-2}%
\end{array}%
$ \\ \cline{2-9}
& \multicolumn{1}{|l|}{LLDP$45$} & $%
\begin{array}{c}
10\times \text{Crude} \\
100\times \text{Mild} \\
10\times \text{Refined}%
\end{array}%
$ & $%
\begin{array}{c}
242 \\
567 \\
3783%
\end{array}%
$ & $%
\begin{array}{c}
29 \\
39 \\
107%
\end{array}%
$ & $%
\begin{array}{c}
1627 \\
3637 \\
23341%
\end{array}%
$ & $%
\begin{array}{c}
271 \\
606 \\
3890%
\end{array}%
$ & $%
\begin{array}{c}
0.49 \\
0.26 \\
0.37%
\end{array}%
$ & $%
\begin{array}{c}
1.0\times 10^{2} \\
5.2 \\
1.1\times 10^{-2}%
\end{array}%
$ \\ \hline
rigid & \multicolumn{1}{|l|}{ode45} & $%
\begin{array}{c}
\text{Crude} \\
\text{Mild} \\
\text{Refined}%
\end{array}%
$ & $%
\begin{array}{c}
19 \\
66 \\
256%
\end{array}%
$ & $%
\begin{array}{c}
2 \\
4 \\
1%
\end{array}%
$ & $%
\begin{array}{c}
127 \\
421 \\
1543%
\end{array}%
$ & $%
\begin{array}{c}
0 \\
0 \\
0%
\end{array}%
$ & $%
\begin{array}{c}
1 \\
1 \\
1%
\end{array}%
$ & $%
\begin{array}{c}
2.7\times 10^{-2} \\
7.5\times 10^{-5} \\
2.0\times 10^{-7}%
\end{array}%
$ \\ \cline{2-9}
& \multicolumn{1}{|l|}{LLDP$45$} & $%
\begin{array}{c}
30\times \text{Crude} \\
22\times \text{Mild} \\
19.7\times \text{Refined}%
\end{array}%
$ & $%
\begin{array}{c}
15 \\
30 \\
111%
\end{array}%
$ & $%
\begin{array}{c}
0 \\
2 \\
1%
\end{array}%
$ & $%
\begin{array}{c}
91 \\
193 \\
673%
\end{array}%
$ & $%
\begin{array}{c}
15 \\
32 \\
112%
\end{array}%
$ & $%
\begin{array}{c}
1.14 \\
0.93 \\
0.83%
\end{array}%
$ & $%
\begin{array}{c}
1.8\times 10^{-2} \\
3.3\times 10^{-5} \\
1.3\times 10^{-7}%
\end{array}%
$ \\ \hline
chm & \multicolumn{1}{|l|}{ode45} & $%
\begin{array}{c}
\text{Crude} \\
\text{Mild} \\
\text{Refined}%
\end{array}%
$ & $%
\begin{array}{c}
679 \\
723 \\
1521%
\end{array}%
$ & $%
\begin{array}{c}
47 \\
16 \\
1%
\end{array}%
$ & $%
\begin{array}{c}
4357 \\
4435 \\
9133%
\end{array}%
$ & $%
\begin{array}{c}
0 \\
0 \\
0%
\end{array}%
$ & $%
\begin{array}{c}
1 \\
1 \\
1%
\end{array}%
$ & $%
\begin{array}{c}
1.1\times 10^{-3} \\
1.1\times 10^{-6} \\
1.3\times 10^{-8}%
\end{array}%
$ \\ \cline{2-9}
& \multicolumn{1}{|l|}{LLDP$45$} & $%
\begin{array}{c}
1.5\times \text{Crude} \\
1.4\times \text{Mild} \\
2.0\times \text{Refined}%
\end{array}%
$ & $%
\begin{array}{c}
146 \\
341 \\
789%
\end{array}%
$ & $%
\begin{array}{c}
1 \\
2 \\
58%
\end{array}%
$ & $%
\begin{array}{c}
883 \\
2059 \\
5083%
\end{array}%
$ & $%
\begin{array}{c}
147 \\
343 \\
847%
\end{array}%
$ & $%
\begin{array}{c}
0.41 \\
1.00 \\
1.10%
\end{array}%
$ & $%
\begin{array}{c}
1.1\times 10^{-3} \\
1.0\times 10^{-6} \\
1.0\times 10^{-8}%
\end{array}%
$ \\ \hline
bruss & \multicolumn{1}{|l|}{ode45} & $%
\begin{array}{c}
\text{Crude} \\
\text{Mild} \\
\text{Refined}%
\end{array}%
$ & $%
\begin{array}{c}
46 \\
148 \\
558%
\end{array}%
$ & $%
\begin{array}{c}
12 \\
13 \\
4%
\end{array}%
$ & $%
\begin{array}{c}
349 \\
967 \\
3373%
\end{array}%
$ & $%
\begin{array}{c}
0 \\
0 \\
0%
\end{array}%
$ & $%
\begin{array}{c}
1 \\
1 \\
1%
\end{array}%
$ & $%
\begin{array}{c}
7.7\times 10^{-2} \\
8.7\times 10^{-6} \\
1.5\times 10^{-8}%
\end{array}%
$ \\ \cline{2-9}
& \multicolumn{1}{|l|}{LLDP$45$} & $%
\begin{array}{c}
4.6\times \text{Crude} \\
1.2\times \text{Mild} \\
3.5\times \text{Refined}%
\end{array}%
$ & $%
\begin{array}{c}
28 \\
101 \\
309%
\end{array}%
$ & $%
\begin{array}{c}
7 \\
13 \\
8%
\end{array}%
$ & $%
\begin{array}{c}
211 \\
685 \\
1903%
\end{array}%
$ & $%
\begin{array}{c}
35 \\
114 \\
330%
\end{array}%
$ & $%
\begin{array}{c}
1.16 \\
1.27 \\
1.05%
\end{array}%
$ & $%
\begin{array}{c}
7.5\times 10^{-2} \\
8.1\times 10^{-6} \\
1.3\times 10^{-8}%
\end{array}%
$ \\ \hline
vdp1 & \multicolumn{1}{|l|}{ode45} & $%
\begin{array}{c}
\text{Crude} \\
\text{Mild} \\
\text{Refined}%
\end{array}%
$ & $%
\begin{array}{c}
59 \\
204 \\
785%
\end{array}%
$ & $%
\begin{array}{c}
10 \\
32 \\
19%
\end{array}%
$ & $%
\begin{array}{c}
415 \\
1417 \\
4825%
\end{array}%
$ & $%
\begin{array}{c}
0 \\
0 \\
0%
\end{array}%
$ & $%
\begin{array}{c}
1 \\
1 \\
1%
\end{array}%
$ & $%
\begin{array}{c}
2.0 \\
2.8\times 10^{-4} \\
5.6\times 10^{-7}%
\end{array}%
$ \\ \cline{2-9}
& \multicolumn{1}{|l|}{LLDP$45$} & $%
\begin{array}{c}
3.2\times \text{Crude} \\
3.2\times \text{Mild} \\
4.0\times \text{Refined}%
\end{array}%
$ & $%
\begin{array}{c}
42 \\
128 \\
461%
\end{array}%
$ & $%
\begin{array}{c}
10 \\
26 \\
10%
\end{array}%
$ & $%
\begin{array}{c}
313 \\
925 \\
2827%
\end{array}%
$ & $%
\begin{array}{c}
52 \\
154 \\
471%
\end{array}%
$ & $%
\begin{array}{c}
1.21 \\
1.12 \\
0.95%
\end{array}%
$ & $%
\begin{array}{c}
1.47 \\
2.4\times 10^{-4} \\
4.1\times 10^{-7}%
\end{array}%
$ \\ \hline
vdp100 & \multicolumn{1}{|l|}{ode45} & $%
\begin{array}{c}
\text{Crude} \\
\text{Mild} \\
\text{Refined}%
\end{array}%
$ & $%
\begin{array}{c}
16916 \\
17516 \\
31253%
\end{array}%
$ & $%
\begin{array}{c}
1074 \\
1540 \\
9%
\end{array}%
$ & $%
\begin{array}{c}
107941 \\
114337 \\
187573%
\end{array}%
$ & $%
\begin{array}{c}
0 \\
0 \\
0%
\end{array}%
$ & $%
\begin{array}{c}
1 \\
1 \\
1%
\end{array}%
$ & $%
\begin{array}{c}
1.9\times 10^{4} \\
0.42 \\
1.2\times 10^{-3}%
\end{array}%
$ \\ \cline{2-9}
& \multicolumn{1}{|l|}{LLDP$45$} & $%
\begin{array}{c}
10\times \text{Crude} \\
39\times \text{Mild} \\
98.5\times \text{Refined}%
\end{array}%
$ & $%
\begin{array}{c}
3780 \\
5026 \\
10719%
\end{array}%
$ & $%
\begin{array}{c}
24 \\
731 \\
37%
\end{array}%
$ & $%
\begin{array}{c}
22825 \\
34543 \\
64537%
\end{array}%
$ & $%
\begin{array}{c}
3804 \\
5757 \\
10756%
\end{array}%
$ & $%
\begin{array}{c}
0.34 \\
0.50 \\
0.53%
\end{array}%
$ & $%
\begin{array}{c}
1.9\times 10^{4} \\
0.26 \\
1.1\times 10^{-3}%
\end{array}%
$ \\ \hline
\end{tabular}%
$} \caption{Code performance in the integration of the nonlinear
examples with similar accuracy.} \label{Table VI ODE}
\end{table}

\section{Discussion}

The results of the previous section show the following: 1) on the
same time partition (Tables \ref{Table I ODE} and \ref{Table II
ODE}), the embedded LLRK formulas are significantly much accurate
than the classical embedded RK formulas of Dormand \& Prince. 2)
with identical tolerances and adaptive strategy (Tables \ref{Table
III ODE} and \ref{Table IV ODE}), the LLDP$45$ code is more accurate
than the ode$45$ code and requires much less time steps for
integrating the whole intervals. For highly oscillatory, stiff
linear, stiff semilinear and mildly stiff nonlinear problems the
overall time of the adaptive LLDP$45$ code is lower than that of the
ode$45$ code, whereas it is similar or bigger for equations with
smooth solution; 3) for reaching similar - but always lower -
accuracy (Tables \ref{Table V ODE} and \ref{Table VI ODE}), the
LLDP$45$ code also requires much less time steps than the ode$45$
code for integrating the whole intervals. In this situation, the
overall time of the adaptive LLDP$45$ code is again much lower than
that of the ode$45$ code for highly oscillatory, stiff linear, stiff
semilinear and mildly stiff nonlinear problems, whereas it is
slightly bigger only for two equations with smooth solution (bruss
and vdp1 examples); and 4) the accuracy of the dense output of the
LLDP$45$ code is, in general, higher than the accuracy of the
ode$45$ code (Tables \ref{Table VII ODE} and \ref{Table VIII ODE}).
\newline
\newline

\begin{table}[H] \centering%
{\scriptsize $%
\begin{tabular}{|c|c|c|c|c|c|}
\hline
Example & Code & Tol & $%
\begin{array}{c}
\text{Time} \\
\text{steps}%
\end{array}%
$ & $%
\begin{array}{c}
\text{Dense} \\
\text{OutPut}%
\end{array}%
$ & $%
\begin{array}{c}
\text{Relative} \\
\text{Error}%
\end{array}%
$ \\ \hline PerLin & \multicolumn{1}{|l|}{ode45} & $\begin{array}{c}
  \text{Crude} \\
  \text{Mild} \\
  \text{Refined}
\end{array}$ & $\begin{array}{c}
  147 \\
  598 \\
  2394
\end{array}$ & $\begin{array}{c}
  589 \\
  2393 \\
  9577
\end{array}$ & $\begin{array}{c}
  10.2\\
  2.6\times 10^{-2} \\
  1.6\times 10^{-3}
\end{array}$
\\ \cline{2-6}
& \multicolumn{1}{|l|}{LLDP$45$} & $\begin{array}{c}
  \text{Crude} \\
  \text{Mild} \\
  \text{Refined}
\end{array}$  & $\begin{array}{c}
  14 \\
  14 \\
  15
\end{array}$ & $\begin{array}{c}
  57 \\
  57 \\
  61
\end{array}$ & $\begin{array}{c}
  2.0\times 10^{-9} \\
  3.0\times 10^{-9} \\
  4.1\times 10^{-9}
\end{array}$ \\
\hline PerNoLin & \multicolumn{1}{|l|}{ode45} & $\begin{array}{c}
  \text{Crude} \\
  \text{Mild} \\
  \text{Refined}
\end{array}$ & $\begin{array}{c}
  105 \\
  411 \\
  1634
\end{array}$ & $\begin{array}{c}
  313 \\
  1169 \\
  4625
\end{array}$ &
 $\begin{array}{c}
  4.8\times 10^{-3} \\
  3.0\times 10^{-6} \\
  2.7\times 10^{-9}
\end{array}$
\\ \cline{2-6}
& \multicolumn{1}{|l|}{LLDP$45$} & $\begin{array}{c}
  \text{Crude} \\
  \text{Mild} \\
  \text{Refined}
\end{array}$ & $\begin{array}{c}
  42 \\
  137 \\
  534
\end{array}$ & $\begin{array}{c}
  101 \\
  293 \\
  1073
\end{array}$ & $\begin{array}{c}
  1.5\times 10^{-3} \\
  8.7\times 10^{-7} \\
  9.2\times 10^{-10}
\end{array}$  \\
\hline StiffLin & \multicolumn{1}{|l|}{ode45} & $\begin{array}{c}
  \text{Crude} \\
  \text{Mild} \\
  \text{Refined}
\end{array}$ & $\begin{array}{c}
  60 \\
  78 \\
  172
\end{array}$ & $\begin{array}{c}
  241 \\
  313 \\
  689
\end{array}$ &
$\begin{array}{c}
  1.1\times 10^{-3} \\
  1.1\times 10^{-6} \\
  8.1\times 10^{-10}
\end{array}$
\\ \cline{2-6}
& \multicolumn{1}{|l|}{LLDP$45$} & $\begin{array}{c}
  \text{Crude} \\
  \text{Mild} \\
  \text{Refined}
\end{array}$ & $\begin{array}{c}
  14 \\
  14 \\
  15
\end{array}$ & $\begin{array}{c}
  53 \\
  57 \\
  61
\end{array}$ & $\begin{array}{c}
  2.7\times 10^{-12} \\
  2.7\times 10^{-12} \\
  2.7\times 10^{-12}
\end{array}$  \\
\hline StiffNoLin & \multicolumn{1}{|l|}{ode45} & $\begin{array}{c}
  \text{Crude} \\
  \text{Mild} \\
  \text{Refined}
\end{array}$ & $\begin{array}{c}
  104 \\
  133 \\
  294
\end{array}$  & $\begin{array}{c}
  417 \\
  533 \\
  1177
\end{array}$   & $\begin{array}{c}
  1.4\times 10^{-2} \\
  3.0\times 10^{-5} \\
  2.6\times 10^{-8}
\end{array}$
 \\ \cline{2-6} & \multicolumn{1}{|l|}{LLDP$45$}
& $\begin{array}{c}
  \text{Crude} \\
  \text{Mild} \\
  \text{Refined}
\end{array}$ & $\begin{array}{c}
  21 \\
  43 \\
  132
\end{array}$  & $\begin{array}{c}
  85 \\
  173 \\
  525
\end{array}$ & $\begin{array}{c}
  6.4\times 10^{-3} \\
  2.9\times 10^{-5} \\
  7.3\times 10^{-8}
\end{array}$  \\
\hline
\end{tabular}%
$}%
\caption{Relative error of the continuous formulas of the codes over
their
dense output after integrating the semilinear examples.}\label{Table VII ODE}%
\end{table}%

\begin{table}[H] \centering%
{\scriptsize $%
\begin{tabular}{|c|c|c|c|c|c|}
\hline
Example & Code & Tol & $%
\begin{array}{c}
\text{Time} \\
\text{steps}%
\end{array}%
$ & $%
\begin{array}{c}
\text{Dense} \\
\text{OutPut}%
\end{array}%
$ & $%
\begin{array}{c}
\text{Relative} \\
\text{Error}%
\end{array}%
$ \\ \hline fpu & \multicolumn{1}{|l|}{ode45} & $\begin{array}{c}
  \text{Crude} \\
  \text{Mild} \\
  \text{Refined}
\end{array}$ &
$\begin{array}{c}
  964 \\
  4474 \\
  19190
\end{array}$ & $\begin{array}{c}
  3857 \\
  17897 \\
  76761
\end{array}$ & $\begin{array}{c}
  9.5\times 10^{2} \\
  19.0 \\
  0.86
\end{array}$ \\ \cline{2-6} &
\multicolumn{1}{|l|}{LLDP$45$} & $\begin{array}{c}
  \text{Crude} \\
  \text{Mild} \\
  \text{Refined}
\end{array}$ & $\begin{array}{c}
  377 \\
  1496 \\
  6021
\end{array}$ & $\begin{array}{c}
  1497 \\
  5985 \\
  24085
\end{array}$ & $\begin{array}{c}
  33.8 \\
  2.8\times 10^{-2} \\
  0.15
\end{array}$
 \\ \hline rigid & \multicolumn{1}{|l|}{ode45} &
$\begin{array}{c}
  \text{Crude} \\
  \text{Mild} \\
  \text{Refined}
\end{array}$ & $\begin{array}{c}
  19 \\
  66 \\
  256
\end{array}$ & $\begin{array}{c}
  77 \\
  265 \\
  1025
\end{array}$ & $\begin{array}{c}
  0.31 \\
  3.4\times 10^{-4} \\
  1.1\times 10^{-6}
\end{array}$  \\
\cline{2-6} & \multicolumn{1}{|l|}{LLDP$45$} & $\begin{array}{c}
  \text{Crude} \\
  \text{Mild} \\
  \text{Refined}
\end{array}$ & $\begin{array}{c}
  16 \\
  53 \\
  201
\end{array}$ & $\begin{array}{c}
  65 \\
  213 \\
  805
\end{array}$ & $\begin{array}{c}
  0.19 \\
  1.7\times 10^{-4} \\
  2.3\times 10^{-7}
\end{array}$  \\
\hline chm & \multicolumn{1}{|l|}{ode45} & $\begin{array}{c}
  \text{Crude} \\
  \text{Mild} \\
  \text{Refined}
\end{array}$ & $\begin{array}{c}
  679 \\
  723 \\
  1521
\end{array}$ & $\begin{array}{c}
  2717 \\
  2893 \\
  6085
\end{array}$ & $\begin{array}{c}
  1.1\times 10^{-3} \\
  1.1\times 10^{-6} \\
  5.7\times 10^{-8}
\end{array}$ \\
\cline{2-6} & \multicolumn{1}{|l|}{LLDP$45$} & $\begin{array}{c}
  \text{Crude} \\
  \text{Mild} \\
  \text{Refined}
\end{array}$ & $\begin{array}{c}
  152 \\
  357 \\
  859
\end{array}$ & $\begin{array}{c}
  609 \\
  1429 \\
  3409
\end{array}$ & $\begin{array}{c}
  9.4\times 10^{-4} \\
  9.2\times 10^{-7} \\
  5.8\times 10^{-8}
\end{array}$ \\
\hline bruss & \multicolumn{1}{|l|}{ode45} & $\begin{array}{c}
  \text{Crude} \\
  \text{Mild} \\
  \text{Refined}
\end{array}$ & $\begin{array}{c}
  46 \\
  148 \\
  558
\end{array}$ & $\begin{array}{c}
  185 \\
  593 \\
  2233
\end{array}$ & $\begin{array}{c}
  8.8\times 10^{-2} \\
  1.0\times 10^{-5} \\
  1.7\times 10^{-8}
\end{array}$ \\
\cline{2-6} & \multicolumn{1}{|l|}{LLDP$45$} & $\begin{array}{c}
  \text{Crude} \\
  \text{Mild} \\
  \text{Refined}
\end{array}$ & $\begin{array}{c}
  36 \\
  105 \\
  396
\end{array}$ & $\begin{array}{c}
  145 \\
  421 \\
  1585
\end{array}$ & $\begin{array}{c}
  6.2\times 10^{-3} \\
  2.4\times 10^{-5} \\
  1.1\times 10^{-8}
\end{array}$  \\
\hline vdp1 & \multicolumn{1}{|l|}{ode45} & $\begin{array}{c}
  \text{Crude} \\
  \text{Mild} \\
  \text{Refined}
\end{array}$ & $\begin{array}{c}
  59 \\
  204 \\
  785
\end{array}$ & $\begin{array}{c}
  237 \\
  817 \\
  3141
\end{array}$ & $\begin{array}{c}
  2.9\times 10^{2} \\
  6.9\times 10^{-4} \\
  4.3\times 10^{-6}
\end{array}$ \\
\cline{2-6} & \multicolumn{1}{|l|}{LLDP$45$} & $\begin{array}{c}
  \text{Crude} \\
  \text{Mild} \\
  \text{Refined}
\end{array}$ & $\begin{array}{c}
  44 \\
  162 \\
  609
\end{array}$ & $\begin{array}{c}
  177 \\
  649 \\
  2437
\end{array}$ & $\begin{array}{c}
  2.25 \\
  2.3\times 10^{-4} \\
  1.9\times 10^{-7}
\end{array}$  \\
\hline vdp100 & \multicolumn{1}{|l|}{ode45} & $\begin{array}{c}
  \text{Crude} \\
  \text{Mild} \\
  \text{Refined}
\end{array}$ & $\begin{array}{c}
  16916 \\
  17516 \\
  31253
\end{array}$ & $\begin{array}{c}
  67665 \\
  70065 \\
  125013
\end{array}$ & $\begin{array}{c}
  2.0\times 10^{4} \\
  0.47 \\
  4.4\times 10^{-3}
\end{array}$ \\
\cline{2-6} & \multicolumn{1}{|l|}{LLDP$45$} & $\begin{array}{c}
  \text{Crude} \\
  \text{Mild} \\
  \text{Refined}
\end{array}$ & $\begin{array}{c}
  3866 \\
  7893 \\
  19887
\end{array}$ & $\begin{array}{c}
  15457 \\
  31573 \\
  79509
\end{array}$ & $\begin{array}{c}
  2.0\times 10^{4} \\
  4.1\times 10^{-2} \\
  2.1\times 10^{-3}
\end{array}$  \\
\hline
\end{tabular}%
$}%
\caption{Relative error of the continuous formulas of the codes over
their
dense output after integrating the nonlinear examples.}\label{Table VIII ODE}%
\end{table}

These simulations results clearly shown that, in the ten examples,
the local linearization of the embedded Runge-Kutta formulas of
Dormand and Prince produces a significant improvement of the
accuracy of the classical formulas. However, this is clearly not a
result that could be expected according to the local truncation
errors given in Theorem \ref{Main Theorem}. This indicates that,
most likely, sharper error estimates could be obtained for the
locally linearized formulas, which is certainly an important open
problem to solve.

Note that the significantly better accuracy of the locally
linearized formulas implies a substantial reduction of the number of
time steps and, consequently, a sensitive reduction of the overall
computation cost in eight of the ten test equations (see Tables
\ref{Table V ODE} and \ref{Table VI ODE}). This indicates that, for
various classes of equations, the additional computational cost of
computing the exponential of a Jacobian matrix at each time step is
compensated for the gain of accuracy. This result certainly agrees
with previous reports in the same direction as that given in
\cite{Shampine97}: "supplying a function for evaluating the Jacobian
can be quite advantageous, both with respect to reliability and
cost".

Further, note that three of the eight test equations for which the
application of locally linearized formulas yields a sensitive
reduction of the overall computation cost are systems of twelve
equations. This illustrates the usefulness of these integrators for
low dimensional problems in general. However, because the locally
linearized formulas (\ref{embedded LLRK45}) are expressed in terms
of the Padé algorithm for computing exponential matrices, it is
expected that they are unable to integrate moderately large system
of ODE with a rational computational cost. In this case, because of
the flexibility in the numerical implementation of the LLRK methods
mentioned in the introduction, the local linearization of the
embedded Runge Kutta formulas of Dormand and Prince can be easily
formulated in terms of the Krylov-type methods for exponential
matrices. In effect, this can be done just by replacing the Padé
formula in (\ref{embedded LLRK45b}) and (\ref{continuousLLRK$45$b})
by the Krylov-Padé formula as performed in \cite{de la Cruz
11,Jimenez09 BIT,Jimenez12 BIT} for the local linearizations schemes
for ordinary, random and stochastic differential equations. In this
way, the Locally Linearized formulas of Dormand and Prince could be
applied to high dimensional ODEs with a reasonable computational
cost \cite{Sotolongo 14}.

On the other hand, we recall that, in order to study the effect of
the local linearization on the conventional RK scheme of Dormand and
Prince, the LLDP$45$ code considered in this work is an exact copy
of the code ode$45$ with the exception of the program lines
corresponding to the embedded and continuous formulas. In this way,
the LLDP$45$ code does not include a number of convenient
modifications that might improve its performance. Some of they are
the following:

\begin{itemize}
\item the initial $h$ at $t_{0}$, which can be estimated by means the exact
second derivative of the solution $\mathbf{x}$ with no extra cost
(as in \cite{Sotolongo 11});

\item online smoothness and stiffness control for estimating the new $h$ at
each step (as, e.g., in \cite{Hairer-Wanner96});

\item the automatic detection of constant Jacobian matrix (as in \cite
{DEUFLHARD92,STEIHAUG79,Zedan90});

\item option for using exact, numerical or automatic Jacobian matrices (as
in \cite{Shampine97, Shampine07, Bischof03});

\item faster algoritms to compute the Pad\'{e} approximation to exponential
matrix (as, e.g., in \cite{Higham05})

\item a parallel implementation of matrix multiplications involved in the exponential matrix evaluations
for taking advantage of the multi core technology available in the
current microprocessors;

\item increase the number of times of the dense outputs: a) up to twelve per
each pair of consecutive times of the partition $(t)_{h}$ with no
extra computation of exponential matrices; or b) up to ninety with
some few extra matrix multiplications;

\item a new continuous formula that replace the current one based on the
continuous RK formula by other based on a polynomial interpolation
of the LLRK formula itself (i.e, derived from the standard way of
constructing continuous RK formulas as in \cite{Hairer-Wanner93});
and

\item change of $h_{max}$, which seems to be too short for semilinear
equations.
\end{itemize}

\section{Conclusions}

In this paper, embedded Locally Linearized Runge-Kutta formulas for
initial value problems were introduced and their performance
analyzed by means of exhaustive numerical simulations. In this way,
the effect that produces the local linearization of the classical
embedded Runge-Kutta formulas of Dormand and Prince were studied. It
was shown that, for a variety of well-known physical equations
usually taken in simulations studies as test equations, the local
linearization of the embedded Runge-Kutta formulas of Dormand \&
Prince produces a significant improvement of the accuracy of
classical formulas, which implies a substantial reduction of the
number of time steps and, consequently, a sensitive reduction of the
overall computation cost of their adaptive implementation.

\section*{Acknowledgment}

The first author thanks to Prof. A. Yoshimoto for his invitation to
the Institute of Statistical Mathematics, Japan, where the
manuscript and its revised version were completed.

\bibliographystyle{elsart-num-sort}



\end{document}